\documentclass[preprint,12pt]{elsarticleA}

\usepackage{amssymb}

\usepackage[ansinew,latin1]{inputenc}
\usepackage[english]{babel}
\usepackage{amsfonts}      
\usepackage{amsmath}
\usepackage{latexsym}
\usepackage{verbatim}    

 \usepackage[caption=false]{subfig}  
   
\def\BARIT#1{{\bar {#1}}}
\def\BARH{\BARIT{H}}
\def\BARHZ{\BARIT{H}^{(0)}}
\def\BARM{\BARIT{\mu}}
\def\BARMZ{\BARIT{\mu}^{(0)}}

\def\LATT{{\Lambda}_N}

\def\LATTC{{\bar{\Lambda}_{M}}}
 
\def\COP{\mathbf{T}}

 \def\AVG#1{{<\!\!#1\!\!>}}
\def\VAR{\mathrm{Var}\,}

\def\Cspace{\{\sigma\in \Sigma\!: \COP\sigma=\eta\}}

\def\BIGO{\mathcal{O}}
\def\VIZ#1{(\ref{#1})}

\def\Aa{\mathcal{A}}
\def\Bb{\mathcal{B}}

\def\Ee{\mathcal{E}}

\def\Rr{\mathcal{R}}

\def\Kk{\mathcal{K}}
\def\Oo{\mathcal{O}}

\def\R{\mathbb{R}}

\def\EXP#1{e^{#1}}

\def\EXPEC#1{{\mathbb{E}\left[{#1}\right]}}

\def\COMMA{\,,}             
\def\PERIOD{\,.}            

\def\VIZ#1{(\ref{#1})}      

\def\BIGO{\Oo}

\newtheorem{algorithm}{Algorithm}
\newtheorem{thm}{Theorem}
\newtheorem{lem}[thm]{Lemma}

\newtheorem{defin}{Definition}

\newproof{pf}{Proof}
\newproof{pot}{Proof of Theorem}
\newdefinition{rmk}{Remark}

\journal{Journal of Computational Physics}

\begin{document}

\begin{frontmatter}

%
%

\title{Multilevel coarse graining and nano--pattern discovery in
many particle stochastic systems}
  
\author[udel]{Evangelia Kalligiannaki\corref{cor1}}
\ead{ekalligi@math.udel.edu}
\author[umass,crete]{Markos Katsoulakis} 
\ead{markos@math.umass.edu}
\author[udel]{Petr Plech\'a\v{c}\corref{cor2}} 
\ead{plechac@math.udel.edu}
\author[udelc]{Dion Vlachos} 
\ead{vlachos@che.udel.edu}

\address[udel]{Department of Mathematical Sciences,
University of Delaware,
Newark, Delaware, United States }
\address[umass]{Department of Mathematics and Statistics, University
of Massachusetts at Amherst, United States}
\address[crete]{ Department of Applied Mathematics, 
University of Crete and Foundation of Research and Technology-Hellas,
Greece }
\address[udelc]{Department of Chemical Engineering,
University of Delaware,
Newark, Delaware, United States }

\cortext[cor1]{Corresponding author}
\cortext[cor2]{Principal corresponding author}

 
%
%

\begin{abstract}
 
In this work we propose a hierarchy of Monte Carlo methods for sampling equilibrium
properties of stochastic lattice systems with competing short and long range
interactions. Each Monte Carlo step is composed by two or more  sub - steps efficiently coupling coarse and microscopic state spaces. The method can be designed to sample the exact or
controlled-error approximations of the target distribution, providing information on
 levels of different resolutions, as well as at the microscopic level. In
both strategies the method achieves significant reduction of the computational
cost compared to conventional Markov Chain Monte Carlo methods. Applications in
phase transition and pattern formation problems confirm the efficiency of the proposed methods.
\end{abstract}


\end{frontmatter}

%
%

\section{Introduction}\label{intro}
Our primary goal in this work is to develop a systematic mathematical and
computational strategy for accelerating microscopic simulation methods with competing short and long range interactions, arising in numerous physical systems for instance in micromagnetics, models of epitaxial growth,  etc.
 We propose the Multilevel Coarse Graining Monte Carlo (Multilevel CGMC)
method, based on a hybrid statistical mechanics and statistics approach. The method introduces  a {\it hierarchy of Markov Chain Monte Carlo methods coupling scales and types of
interactions} that can sample the exact or controlled error approximations of 
Gibbs measures $\mu_{N,\beta}(d\sigma) = Z_{N}^{-1}\EXP{-\beta H_N(\sigma)}\,P_N(d\sigma)$  defined on a high dimensional space, $ N>>1$, that  can be easily generalized to any probability 
measure with similar properties.
  It is a method of constructing efficient {\it proposal measures}  in 
Metropolis sampling using coarse--graining techniques, see also \VIZ{muzero} below,
aiming at reducing the rejection rate and the computational complexity. The key
idea is a decomposition of the sampling distribution to a product measure \begin{equation}\label{Decompose}
\mu_{N,\beta}(d\sigma) = \BARMZ_{M,\beta}(d\eta) \nu_r(d\sigma|\eta)\COMMA
\end{equation}
with $ \eta := \COP \sigma $ a variable with less degrees of freedom compared to $\sigma$, defined by a projection map $\COP:\Sigma_N \to \BARIT{\Sigma}_M,\ M<N$. $\BARMZ_{M,\beta}(d\eta)$ is a measure with a simple explicit representation approximating the marginal  $\BARM_{M,\beta}(d\eta) = \mu_{N,\beta} \circ \COP^{-1}(d\eta)$ and $\nu_r(d\sigma|\eta)$ is a uniquely defined (prior) measure, responsible for {\it reconstructing } variables $\sigma $ given $ \eta$, \cite{KPR}. Such a two--level measure decomposition can be trivially extended to a multi--level setting where \VIZ{Decompose} can include different resolution
levels interpolating between a coarser level and the microscopic one $\sigma$.

 In the following we describe a  Monte Carlo step of a two--level  CGMC method:
\begin{enumerate}
\item Sample $\eta$ from  $\BARMZ_{M,\beta}(d\eta)$, using Coarse Grained (CGMC)
samplers \cite{KMV, KMV1}.  Appropriate coarse grained measures have been evaluated in
earlier work via cluster expansions that can be easily constructed  with available analytical
error estimates ensuring that such approximations are controllable \cite{KPRT}. 
\item Conditioned on such $\eta$, obtained in Step 1, we sample $\nu_r(d\sigma|\eta)$ via an
Accept/Reject method.
\end{enumerate}
A schematic description of this procedure can be seen in Figure \ref{fig:0}.
Better proposals constructed in Step 1 will lead to fewer rejections in Step 2,
furthermore, there is no need to consider all possible microscopic proposals
since at the coarse step we do a first
screening. Instead of further approximating the sampling measure, as is done with cluster expanding \VIZ{rg} keeping the typically computationally expensive multi-body higher order terms \cite{ KPRT3}, we use the hybrid statistical and statistics approach that
construct $\nu_r(d\sigma|\eta)$. Even when CGMC  provide less accurate
approximations, the multilevel CGMC approach can refine the results by the
Accept/Reject step in the finer space. A key issue is how to sample conditionally from $\nu_r$ in an efficient manner in Step 2, picking an appropriate coarse proposal in Step 1.

The necessary ingredient for the applicability of the method
is a decomposition of the form \VIZ{Decompose}, which includes a possibly less
accurate coarse--grained measure and the correcting accept/reject Step 2 above.
This formulation can make the proposed method extentable to off--lattice systems where various coarse--graining schemes are already available 
\cite{kremerplathe,vagelis, AKPR}, although without controlled--error
approximations.
In such off--lattice systems we typically have two main features: a presence of
short and longer interactions, as well as comparable energy and entropy, hence
fluctuations are expected to be important in the modelling and simulation.

Systems with smooth long or intermediate range interactions are well approximated by coarse graining techniques \cite{KMV1, KPR, KT}, and  CGMC are reliable simulation methods with controlled--error approximations,
both for observables and loss of information \cite{KPS,
KPRT}. Furthermore, models where only short--range interactions appear are
inexpensive to simulate with conventional methods.
However, when both short and long range interactions are present, the conventional methods  become  prohibitively expensive, and  CG error estimates are not applicable.
The proposed method can handle such systems efficiently by either    compressing only the long range interactions, and sample with CGMC  with low computational cost, at  Step 1 and  include the  short range part at the accept/reject Step 2,  or compress   all types of interactions for Step 1, and correct   appropriately  in   Step 2.

A wide literature exists on sophisticated Markov Chain Monte Carlo (MCMC) methods designed  to accelerate  simulations for large systems,  applying   for example  parallelising techniques   and/or   constructing  good   first approximations (proposals) in Metropolis sampling \cite{Duane}, \cite{Neal}. In \cite{EHL} Efendiev et.al., propose the Preconditioning MCMC, a two stage  Metropolis  method, applied to inverse problems of subsurface characterization. Our
algorithm shares the same idea   of constructing a proposal density based on meso /macro-scopic properties of the model studied and taking advantage of the first stage rejections.
Several methods  where  the  trial density is built up sequentially with stage-wise rejection decision appear, \cite{CEP},\cite{LIU}. There are also some similarities with
simulated sintering, and transdimensional MCMC, see \cite{LiuSabatti1,LIU} and
references therein. However, the novelty of our method lies on the construction of the variable
dimensionality (and level of coarse--graining) state spaces and the
corresponding Gibbs measures relies on statistical mechanics tools that allow a systematic
control of the error from one level of coarse--graining to the next. The interplay of different levels of compressed spaces appears also in spatial multi--grid methods coupled with CGMC sampling, studied in \cite{sinno} \cite{vlachosmg}, for accelerating lattice kinetic Monte Carlo simulations where however the proposed methods are not exact and do not necessarily provide controlled--error approximations.  Various  attempts  appear on parallelising Monte Carlo simulations,  based on a parallel resolution,  such as parallel kinetic CGMC \cite{AKP} and a
combination of CGMC and parallel tempering \cite{MTP}.
In a follow up work we extend our framework for multilevel CGMC
to accelerate sampling of  evolution processes on large lattice systems \cite{KKP2}. There we develop    multilevel Kinetic Monte Carlo algorithms, based on the known  from coarse graining techniques different level  approximating processes.


In Section \ref{SLS} we introduce microscopic lattice systems and provide a brief review of
coarse graining  methods. We also provide the Metropolis-Hastings MCMC
method for numerical simulations and describe microscopic processes. 
We introduce the multilevel CGMC (ML CGMC) method in Section~\ref{sectionEquil}  and provide the mathematical analysis that    ensures the theoretical validity of the method.    Section~ \ref{complexity}  provides a full comparison of the computational complexity between the classical and the multilevel Metropolis introduced here. 
More specifically, in Theorem~\ref{main} we prove comparative estimates on their mixing times via   spectral gap estimates. In turn such spectral estimates   are obtained through suitable  relative bounds on  their respective Dirichlet forms. Concluding, our  analysis shows that  ML CGMC has substantial savings over the classical MH algorithm, generating cheaply proposals with small or controllable rejection rates.
Sections \ref{applications} and \ref{patterns} give example applications of the
multilevel method in canonical and microcanonical sampling. In Section
\ref{benchmark} a benchmark example is employed in order to demonstrate an
explicit application of Theorem~\ref{main}. In subsection \ref{nontrivial} an order one improvement of the coarse graining error in a phase transition
regime is achieved when applying the multilevel CGMC method with a potential  splitting approach, for a Kac type potential with algebraic decay. Finally in
Section \ref{patterns}  we study nanopattern formation in surface diffusion
induced by a Morse type potential, and verify   that the  proposed method  can provide  correctly microscopic details.

%
%

\section{Stochastic lattice systems at equilibrium}\label{SLS}
We consider an Ising-type system on a periodic $d$-dimensional lattice $\LATT$
with $N=n^d$ lattice sites. At each site $x \in \LATT$ we define a spin
variable $\sigma(x)$ taking values in a finite set. For instance, in a lattice gas model
$\sigma(x)\in\{0,1\}$ describes that the site $x$ is vacant or occupied by an atom. The state of the system is described by a configuration $\sigma\in \Sigma_N = \{0,1\}^{\LATT}$. The interaction energy of the system, e.g., interacting particles in the lattice gas model, is defined by
the Hamiltonian $H_N$. We assume systems where the particles interact only
through a pair-wise potential and thus the Hamiltonian takes the form
$\sigma=\{\sigma (x):x\in\LATT\}$
\begin{equation}\label{GenHam}
H_N(\sigma) = -\frac12\sum_{x \in \LATT}\sum_{y\not= x} J(x-y) \sigma
(x)\sigma(y) +
  \sum_{x\in \LATT}h(x)\sigma(x) \COMMA
\end{equation}
where $h$ is an external field. 
Equilibrium states at the inverse temperature $\beta$ are described by the
(canonical) Gibbs probability measure on the space $\Sigma_N$
\begin{equation}\label{gibbs}
\mu_{N, \beta}(d \sigma) = Z_{N}^{-1}\EXP{-\beta
H_N(\sigma)}\,P_N(d\sigma)\COMMA
\end{equation}
where $Z_{N}$ is the normalizing constant (partition function). Furthermore,
the product Bernoulli distribution $P_N(d\sigma)$ is the {\it prior distribution} on
$\LATT$ representing distribution of states in a non-interacting system, or
equivalently at $\beta=0$, when thermal fluctuations-disorder-associated with the product structure of $P_N(d\sigma)$ dominates. By contrast at zero temperature, $\beta=\infty$,
interactions and hence order, prevail. Finite temperatures, $0<\beta <\infty$,
describe intermediate states, including possible phase transitions between
ordered and disordered states.
 
The coarse-graining techniques have been developed in order to study the
behaviour in the regimes when the size of the system $N\to\infty$. In the series of
papers \cite{KMV1, KMV, KV} the authors initiated the development of {\em
coarse-graining} (CG) as a computational tool for accelerating Monte Carlo
simulations of stochastic lattice dynamics. The coarse-grained model is
constructed on a coarse grid $\LATTC$ by dividing $\LATT$ into $M$ coarse
cells, each of which contains $Q$ (micro-)cells,  typically $Q=q^d$   with the coarse-graining
scale $q$ in each dimension. Each coarse cell is denoted by
$C_k$, $k \in \LATTC$. A typical choice for the coarse variable in the
context of Ising-type models is the block-spin over each coarse cell $C_k$,
\begin{equation*}
\eta:=\left\{ \eta(k)=\sum_{x \in C_k} \sigma(x) \,:\,k \in
\LATTC\right\} \COMMA
\end{equation*}
defining the coarse graining map $ \COP: \Sigma_N \to  \BARIT\Sigma_M, \ \ \COP\sigma=\eta$.
The exact coarse-grained Gibbs measure is given (with a slight abuse
of notation) by
$\BARM_{M,\beta}= \mu_{N,\beta} \circ \COP^{-1}$, written in a
more convenient form
\begin{equation}\label{cg_gibbs}
     \BARM_{M,\beta}(d{\eta})=\frac{1}{\BARIT{Z}_M}
     \EXP{-\beta\BARH_M({\eta})} \BARIT{P}_M( d{\eta})\COMMA
\end{equation}
where $ \BARIT{P}_M( d{\eta}) = P_N \circ \COP^{-1} $ defines the
coarse-grained prior measure. The {\em exact coarse-grained Hamiltonian} is
defined by the {\em renormalization group map}, see, e.g., \cite{Golden},
\begin{equation}\label{rg}
     \EXP{-\beta \BARH_M (\eta)} = \EXPEC{\EXP{-\beta H_N}|{\eta}}
= \int \EXP{-\beta H_N(\sigma)} P_N(d\sigma|{\eta})\PERIOD\end{equation}
The conditional prior $P_N(d\sigma|{\eta})$ is the probability of having a
microscopic configuration $\sigma$, given a coarse configuration $\eta$.

Although typically $\BARIT{P}_M( d{\eta})$ is easy to calculate the exact
computation of the coarse-grained Hamiltonian $\BARH_M(\eta)$ given by \VIZ{rg}
is, in general, impossible even for moderately small values of $N$. Therefore
suitable approximations have to be constructed. An initial approximation
 can be then proposed by the technique of {\em cluster
expansions}, \cite{Simon} providing improved approximation at suitable phase
regimes. The corresponding first order CG Hamiltonian is explicitly given, \cite{KPRT},
\begin{eqnarray}\label{Hzero}
\BARHZ(\eta) =& -&\frac12\sum_{k\in\LATTC}\sum_{l\neq k} \BARIT{J}(k,l)\eta(k)\eta(l) -\frac12 \BARIT{J}(0,0)\sum_{k\in\LATTC}\eta(k)(\eta(k)-1) \nonumber \\ &+&\sum_{k\in\LATTC}\BARIT{h}\eta(k)\COMMA 
\end{eqnarray}
where the coarse-grained interactions are evaluated explicitly by averaging
over the cells $k,l\in\LATTC$
\begin{equation*}
\BARIT{J}(k,l)=\frac{1}{q^2}\sum_{x\in C_k}\sum_{y\in C_l}
J(x-y)\COMMA\;
\BARIT{J}(k,k)=\frac{1}{q(q-1)}\sum_{x\in C_k} \sum_{y\in C_k,y\neq
x}J(x-y)\COMMA
\end{equation*}
defining the coarse grained  Gibbs measure
\begin{equation}\label{muzero}
\BARMZ_{M,\beta}(d\eta) = \frac{1}{\bar Z^{(0)}} \EXP{ -\beta \BARHZ (\eta)}\BARIT{P}_M( d{\eta})\PERIOD
\end{equation}
  The coarse-graining of systems with purely long- or intermediate-range
interactions was studied using cluster expansions in \cite{KPRT, AKPR, KPR,
KPRT2}. In many applications the long-range potentials
exhibit scaling property
 \begin{equation}\label{microJL}
J(x-y)=L^{-d}V\left({\frac{n}{L}|x-y|}\right)\COMMA\;\; x, y \in
\LATT\COMMA
\end{equation}
where $V \in C^1([0,\infty))$ and it is normalized to ensure that the strength
of the potential $J$ is essentially independent of $L$, i.e., $\sum_{x\neq 0}
J(x)\sim \int_0^\infty V(r) dr $. The constant $L$ can be interpreted as a
(characteristic) interaction range of the potential. For example, if we have
$V$ with properties $V(r)=V(-r)$, $V(r)=0$, $|r| > 1$, then a spin at the site $x$
interacts with its neighbours which are at most $L$ lattice points away from
$x$.
One of the results therein is on deriving  error estimates  in terms of the specific relative
entropy $\Rr(\mu| \nu):=N^{-1}\sum_\sigma \log\big\{\mu(\sigma)/
\nu(\sigma) \big\}
\mu(\sigma) \quad 
$ between  the corresponding equilibrium Gibbs measures. Note that the scaling factor $N^{-1}$ is related to the extensivity of the system, hence the proper error quantity that needs to be tracked is the  loss of information {\it per particle}. 
\begin{equation} 
 \Rr(\BARMZ_{M,\beta}|\mu_{N,\beta}\circ \COP^{-1}) = 
    \BIGO(\epsilon^{2})\COMMA \quad\quad \epsilon \,\equiv\,    \beta  \|\nabla V\|_1 \left(\frac{q}{L}\right)\, ,\label{errb}
\end{equation}
 
{\it Systems with short and long range interactions.}
 One of our goals in this work is to study  systems where in addition to the long-range potential we have a short range
 \begin{equation}\label{microK}
K(x-y)=S^{-d}V_s\left({\frac{n}{S}|x-y|}\right)\COMMA\;\; x, y \in \LATT\COMMA
\end{equation}
with $V_s$ having similar properties as $V$ in \VIZ{microJL} and $ S \ll L$ distinguishing the short and long range nature of interactions.  A typical case of a short--range potential is encountered in the nearest-neighbour Ising model where    $K(x-y) = K=constant$, for $|x-y| = S=1$ and zero otherwise.   The new Hamiltonian including the contributing energy from the short range potential $ H_s(\sigma)$ and the long range part $H_l(\sigma) $ is
\begin{equation}\label{HamSplit}
H_N(\sigma) =H_s(\sigma) + H_l(\sigma)\PERIOD
\end{equation}
Study of equilibrium properties of stochastic lattice systems mainly accounts
for the evaluation of averages over the coarse-grained $ \BARM_{M,\beta}(d\eta)$
or the microscopic $\mu_{N,\beta}(d\sigma)$ Gibbs measures of observables
$\phi:\Sigma_N\to\R$, i.e.,
\begin{equation*}
\EXPEC{\phi}= \int_{\Sigma_N} \phi(\sigma)
\mu_{N,\beta}(d\sigma)\PERIOD
\end{equation*}
Numerical methods evaluating equilibrium averages are the Markov Chain Monte Carlo
(MCMC) methods, among which the most widely used is the Metropolis -- Hastings (MH)
method \cite{MRR, HA}. 
Metropolis--Hastings algorithm generates proposals $\sigma'$, for the evolution from the
configuration $ \sigma$ to $\sigma'$, that are defined by  the {\it
proposal} probability transition kernel $\rho(\sigma',\sigma)$. The proposal $\sigma'$ is accepted with probability $\alpha(\sigma,\sigma')$ or rejected with the probability $1-\alpha(\sigma,\sigma')$. Let $X_0=\sigma_0$ be an arbitrary initial configuration, the $n$-th iteration of the algorithm consists of the following steps
\begin{algorithm}[Metropolis-Hastings algorithm]\label{MH}
\noindent Given $X_n=\sigma$  
\begin{description}
\item[Step 1] Generate $X'_n=\sigma' \sim \rho(\sigma',\sigma)$.
\item[Step 2] Accept-Reject
\begin{equation*}
   X_{n+1}= \begin{cases}
X'_n=\sigma'\; &\mbox{with probability $\alpha(\sigma,\sigma')$,} \\
X_n =\sigma \; &\mbox{with probability $1-\alpha(\sigma,\sigma')$,} \\
           \end{cases}
\end{equation*}
with the acceptance probability depending on the energy difference between
configurations $\sigma$ and $\sigma'$, $ \Delta H_N(\sigma,\sigma')=
H_N(\sigma') - H_N(\sigma)$,
\[
\alpha(\sigma,\sigma')=\min\left\{1,\exp\{-\beta \Delta H_N
(\sigma,\sigma')\}
\frac{ \rho(\sigma',\sigma)}{\rho(\sigma,\sigma')} \right\}\PERIOD\]
\end{description}
\end{algorithm}
The algorithm generates an ergodic
Markov chain $\{X_n\}$ in the state space $\Sigma_N$, with the stationary measure
$\mu_{N,\beta}(d\sigma)$. Ergodicity ensures the convergence of empirical
averages $ \frac{1}{n}\sum_{i=1}^n\phi(X_i)$ to the desired mean $ \EXPEC{\phi} $, for any $ \phi\in L^1(\mu_{N,\beta})$.  
It is easy to deduce the    probability transition kernel associated to MH Algorithm~\ref{MH} 
\begin{equation}\label{KC}
\Kk_c(\sigma,\sigma')=\alpha(\sigma,\sigma')\rho(\sigma,\sigma')+
\left[1-\int_{\Sigma}\alpha(\sigma,\sigma')\rho(\sigma,\sigma')d\sigma'\right]
\delta(\sigma'-\sigma)\PERIOD
\end{equation}
where $\delta$ denotes the Dirac measure.
 
Depending on whether one considers microcanonical or canonical ensemble the
ergodic Markov chain $\{X_n\}$ can be defined by spin--exchange dynamics that
preserve the order parameter  or spin--flip dynamics 
respectively \cite{LIGGETT}, \cite{DMP}. In the spin-exchange the  proposed new configuration
$\sigma'=\sigma^{(x,y)}$ is obtained from $\sigma$ by interchanging the spins
at $x$ and $y$, for  nearest-neighbour sites $x$ and $y$.
 \begin{equation*}
   \sigma^{(x,y)}(z)= \begin{cases}
       \sigma(y)\COMMA & \mbox{when $z=x$,}  \\
       \sigma(x)\COMMA & \mbox{when $z=y$,} \\
       \sigma(z)\COMMA & \mbox{otherwise. } \\
    \end{cases}
 \end{equation*}
Analogously for the spin-flip    $\sigma' = \sigma^{(x)}$ is obtained
from $\sigma$ by flipping the spin value at site $x$
 \begin{equation*}
 \sigma^{(x)}(z)= \begin{cases}
                   1-\sigma(x)\COMMA & \mbox{ when $z=x$,}\\
                   \sigma(z)\COMMA & \mbox{ otherwise.} \\
                  \end{cases}
 \end{equation*}

%
%
\section{The multilevel CGMC Metropolis-Hastings method}\label{sectionEquil}
 We present in detail and generalize the Coupled Metropolis--Hastings method originally proposed in \cite{KKP}.  We introduce a multi--level Metropolis--Hastings method, where each level corresponds to a configuration space resolution level, and provide the associated  mathematical analysis  that   ensures the theoretical validity of the method. 
\begin{figure}[ht]
\centering
\includegraphics[angle=-0,width=0.6\textwidth,
height=0.45\textheight]{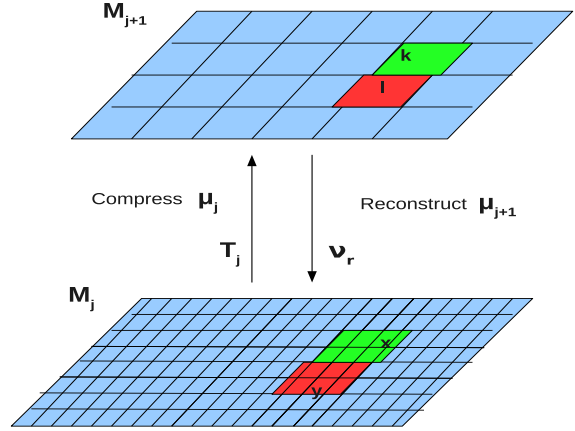}
\caption{  Schematic of a two-level ML CGMC. Information exchange
between coarser and finer resolutions.  Step 1: Sample  $ \mu_{j+1}(d\eta_{j+1})$,   Step 2: Reconstruct  with  $\nu_{r}(d\eta_j|\eta_{j+1})$ such that  $  \mu_{j}(d\eta_j) =  \mu_{j+1}(d\eta_{j+1}) \nu_{r}(d\eta_j|\eta_{j+1})$. }
\label{fig:0}
\end{figure}
 
 
Let $\{\bar{\Sigma}_{M_j}\}_{j=0}^{I}$ denote a hierarchy of coarse spaces
derived from the microscopic space $\Sigma_N=:\BARIT{\Sigma}_{M_0} $ by a family
of mappings $\COP_j: \BARIT{\Sigma}_{M_j} \to \BARIT{\Sigma}_{M_{j+1}}$, $\COP_j \eta_j =
\eta_{j+1}$, for $N=M_0>M_1>\dots>M_I$. The variables $\eta_j$ denote
configurations on spaces $\BARIT{\Sigma}_{M_j} $, with $\eta_0 =: \sigma$ referring
to the microscopic variable on $\Sigma_N$. The method is composed by  a sequence of $ I+1$
Metropolis-Hastings steps each one designed to generate samples from $\bar
\Sigma_{M_j}$ given a coarser sample from $\BARIT{\Sigma}_{M_{j+1}}$. Properly
constructed measures on $\bar{\Sigma}_{M_{j}},\; j>0$ form the basis for constructing
efficient proposal kernels for Metropolis-Hastings algorithm allowing sampling of large systems.
The interplay between  different resolution  spaces  is  controlled by  the hierarchy of projection mappings $\COP_j$ and  corresponding  inverse mapping procedures (reconstruction).  Appropriate reconstruction  measures   $ \nu_{r,j}(d\eta_j|\eta_{j+1})$ are constructed in view of   decompositions  similar to \VIZ{Decomp},  for the each pair of  Gibbs measures $\mu_j(d\eta_j) $ and $ \mu_{j+1}(d\eta_{j+1})$, see Figure \ref{fig:0}.
 For the sake of exposition we describe in
detail the {\it two-level CGMC} method, in which two different scale state
spaces $\BARIT{\Sigma}_{M_{j_1}} = \Sigma_N$ and $\BARIT{\Sigma}_{M_{j_2}}=\Sigma_M$, $M<N$ are involved. 
 
The ML CGMC method is a hybrid statistical mechanics and statistics approach
based on a product measure decomposition of the target distribution
\begin{equation}\label{Decomp}
\mu_{N,\beta}(d\sigma) = \BARMZ_{M,\beta}(d\eta)\nu_r(d\sigma|\eta)\PERIOD
\end{equation}
The principal idea is that computationally inexpensive CG
simulations will reproduce the large scale structure and subsequently
microscopic information will be added through the {\em microscopic
reconstruction}. The two levels of the method  are described by:
\begin{enumerate}
\item {\it Marginal approximation and CGMC sampling.} 
 $\BARMZ_{M,\beta}(d\eta)$, given in \VIZ{muzero},  an  approximation of marginal  $\BARM_{M,\beta} = \mu_{N,\beta}\circ \COP^{-1} $ , is   described explicitly   with a controllable  error   \cite{KPRT}.   
\item {\it Microscopic reconstruction.} Conditioned on $\eta$, obtained by sampling from
$\BARMZ_{M,\beta}(d\eta)$, we sample $\nu_r(d\sigma|\eta)$ via a simple
inverse--mapping distribution $\mu_r(d\sigma|\eta) $ and an Accept/Reject  method. Measure $\nu_r(d\sigma|\eta) $ is responsible for the  {\it reconstruction}
procedure of the fine configuration $\sigma$ constrained on the coarse
configuration $\eta$. 
\end{enumerate}
More specifically, reconstruction is the reverse procedure  to coarse--graining, i.e., reproducing {\it microscopic} properties directly from CG simulations.  A detailed discussion on reconstruction can be found in \cite{KPR, KKP}.   Note that $\nu_r(d\sigma|\eta) $ is a finite measure, uniquely defined by \VIZ{Decomp} given $\mu_{N,\beta}(d\sigma)$ and $\BARMZ_{M,\beta}(d\eta)$,
  
In view of Decomposition \VIZ{Decomp}  we propose the two--level CGMC algorithm  composed of two Metropolis--Hastings steps. The first step samples the measure
$\BARMZ_{M,\beta}(d\eta)$  \VIZ{muzero}, on the coarse space $ \BARIT{\Sigma}_M$, using
an arbitrary proposal transition kernel $\bar{\rho}(\eta,\eta')$ to produce
coarse trial samples $\eta'$. The second step, {\it performed only if the coarse trial
sample is accepted}, consists of the simple reconstruction from the coarse
state $\eta'$ with $ \mu_r(d\sigma|\eta)$ and an accept/reject method.
If a trial coarse sample is rejected the second step is not performed and no
computational time is wasted on checking fine trial samples that are most
likely to be rejected. The two--level CGMC algorithm is defined as follows:
Let $Y_0=\sigma_0$ be an arbitrary initial configuration, for
$n=0,1,2,\dots$,
\begin{algorithm}[Two-level   CGMC MH   Algorithm]\label{mCGMC} 
\noindent 
 Given $Y_n=\sigma$ 
\begin{description}
\item[Step 1]  Compute the coarse configuration $  \eta=\COP\sigma$.
\item[Step 2]  Generate a coarse sample  $\eta'  \sim \bar{\rho}(\eta,\eta')$.
\item[Step 3]  Coarse Level Accept-Reject. \\
	 Accept $\eta'$ with probability: 
	\begin{equation}\label{coarseAcc}
		\alpha_{CG}(\eta,\eta')=\min\left\{1,\EXP{-\beta \Delta \BARHZ(\eta,\eta')} \frac{ \bar{\rho}(\eta',\eta)}{
\bar{\rho}(\eta,\eta')}\right\}\PERIOD
	 \end{equation} 
	  {\bf if} $\eta'$ is accepted then proceed to Step 4, \\
	  {\bf else} generate a new coarse sample Step 2.
\item[Step 4]  Reconstruct $\sigma'$ given  the coarse trial $\eta'$,       
	\[
            \sigma' \sim  \mu_r(\cdot|\eta')\PERIOD
	\] 
\item[Step 5] Fine Level Accept-Reject. \\
	Accept $\sigma'$ with probability
	\begin{equation}\label{fineAcc}
\alpha_f(\sigma,\sigma')=\min\left\{1,\EXP{ -\beta\left[ \Delta
H_N(\sigma, \sigma') - \Delta \bar
H^{(0)}(\eta,\eta')\right]}\frac{\mu_r(\sigma|\eta)}{\mu_r(\sigma'|\eta')}\right\}\PERIOD
	 \end{equation}
\end{description}
\end{algorithm}
where $ \Delta H_N(\sigma,\sigma')= H_N(\sigma') - H_N(\sigma)$, and $\Delta
\BARHZ(\eta,\eta')= \BARHZ(\eta')- \BARHZ(\eta)$.

In terms of Metropolis--Hastings Algorithm~\ref{MH}, 
the method is  generating trial samples $\sigma'$   from  proposal kernel $\rho(\sigma,\sigma')$, with stationary measure $\BARMZ_{M,\beta}(d\eta)\mu_r(d\sigma|\eta)$,  that depends on the statistical mechanics properties of the system.  This fact  leads to an increase of the acceptance rate of the MH method, see for example Figure \ref{fig:2.1a}.
Indeed, consider  the  canonical ensemble for $x\in\LATT$ such that $x\in C_k$, $k\in\LATTC $  with a   potential $ J(x-y)$  \VIZ{microJL},   where the following estimate holds \cite{KMV},
\begin{equation*}
\Delta H_N(\sigma,\sigma^{(x)}) - \Delta \BARHZ(\eta,\eta^{(k)})  = \BIGO\left(Q/ (2L+1)^d\right)\COMMA
\end{equation*}
with $ \eta^{(k)} =\COP\sigma^{(x)}$.  
This estimate shows that the second level acceptance probability $ \alpha_f(\sigma,\sigma')$ is controlled by the coarsening parameter,  i.e. it  is close to $1$   for $Q/(2L+1)^d<<1$ and most of the    coarse samples  entering the second level  will be accepted. Thus, for fixed $L$,
varying the coarsening parameter  we can control the effectiveness of the method as a
balance between acceptance rate and computational cost. That is between small values of
$Q$, leading to high acceptance rate, and larger $ Q$, leading to lower computational
cost at sampling at the coarse space. We elaborate with this issue in detail in
Section~\ref{complexity}.

\subsection{Examples of multi-level decompositions}\label{Decompositions}
In this section we elaborate   on  how to select the prior measure  $\nu_r(d\sigma|\eta)$,  that furthermore determines how to perform a multi--level decomposition of the Gibbs measure.
The selection of   $\nu_r(d\sigma|\eta)$   depends on two features (a) the  choice of the coarse  grained measure $ \BARMZ_{M,\beta}(d\eta)$, that is whether we compress the entire interaction potential or split to a short and long range part  and  (b) on sampling the exact  or an approximation  of the microscopic  measure $ \mu_{N,\beta}(d\sigma)$.

\noindent
{\it  Corrections.}
With this strategy   all types of   interactions  are compressed  and  incorporated into the first level of sampling with  $ \BARHZ(\eta)$,  as is depicted in the algorithmic description of the method  Algorithm~\ref{mCGMC}. The reconstructive and corrective characteristics of  $ \nu_r(d\sigma|\eta)$  appear  explicitly    in the following  formulation, 
that is reconstruct with the uniform conditional distribution    $\mu_r(d\sigma|\eta) = P_N(d\sigma|\eta)$  and correct with $ H_N(\sigma) -\BARHZ(\eta)$,
\begin{equation*}
 \nu_r(d\sigma|\eta) = \BARIT{Z}^{(0)}_{M}Z^{-1}_N\exp\{ -\beta [H_N(\sigma) -\BARHZ(\eta)] \}P_N(d\sigma|\eta)\PERIOD
\end{equation*}

\noindent
{\it Potential splitting.}
An alternative approach is considered  based on  splitting  the inter particle interactions into {\it short}-- and {\it long}--range terms. A decomposition of the coarse-graining of the interaction potential can be justified and optimized by  known error estimates, see \cite{AKPR}. These estimates suggest a natural way to split   the potential into a short--range piece $ K(x-y)$ with possible singularities and a locally integrable (or smooth) long--range decaying component $ J(x-y)$.

Here we suggest to sample on the coarse step  according to  the effective Hamiltonian $\BARHZ_l(\eta)$,  corresponding only to the long--range depended energy $H_l(\sigma)$, as in \VIZ{HamSplit}, that  suggests a  rearrangement of decomposition  \VIZ{Decomp} where  
\begin{equation*}
 \nu_r(d\sigma|\eta) = \BARIT{Z_l}^{(0)}_{M}Z^{-1}_N\exp\{ -\beta [H_s(\sigma) + H_l(\sigma) -\BARHZ_l(\eta) ] \}P_N(d\sigma|\eta)\PERIOD
\end{equation*}
In the  two--level CGMC sampling  the costly long--range part is involved only in the coarse updating where the number of operations to calculate energy differences is reduced and  coarsening of the short-range potential is avoided \cite{KPRT3}.
 
{\it Approximate CG.} In many applications where meso/macroscopic information is sufficient  CGMC sampling is reliable for  long--range potentials and the error when neglecting terms $ \Delta H_l(\sigma,\sigma') - \Delta \bar H_l^{(0)}(\eta,\eta')$ is small.  In this strategy we suggest to neglect these terms, despite the encountering of the approximating error,   benefiting from    a  further reduction of  the computational complexity, see Section~\ref{complexity}.
As a result  Algorithm~\ref{mCGMC} is sampling from a probability measure approximating $\mu_{N,\beta}(d\sigma)$,
\begin{equation*}
\mu^{(0)}_{N,\beta}(d\sigma) \propto \exp\{-\beta H_s(\sigma) - \beta H_l^{(0)}(\eta)\} P_N(d\sigma|\eta)\BARMZ_{M,\beta}(d\eta) \PERIOD
\end{equation*}

\subsection{Reversibility of the Multilevel CGMC } 
Mathematical analysis for  the two--level CGMC method is studied in this section for a broad class of probability measures. The method  can be generalized to  the sampling   of any probability measure $\mu(d\sigma)$  on a countable configuration space $\Sigma$, by coupling properly configurations between a hierarchy of coarser and finer spaces. Since almost any
probability density $\mu(\sigma)$ can be written in the form $ \EXP{-H(\sigma)}$ the method can applied to any model for which one can properly define a function $H(\sigma)$ and the hierarchy of coarse spaces and densities. The following arguments are straightforward but necessary to prove that the algorithm samples processes with the desired stationary
measure.  Algorithm~\ref{mCGMC} is defined by the acceptance probabilities
 \begin{equation}\label{gcoarseAcc}
\alpha_{CG}(\eta,\eta')=\min\left\{1,\frac{\BARMZ(\eta')\bar{\rho}
(\eta',\eta)}{\BARMZ(\eta)\bar{\rho}(\eta,\eta')}\right\}\COMMA
 \end{equation}
 and
 \begin{align}\label{gfineAcc}
\alpha_f(\sigma,\sigma')=\min\left\{1,\frac{\mu(\sigma')
\BARMZ(\eta)\mu_r(\sigma|\eta)}{\mu(\sigma)\BARMZ(\eta')
\mu_r(\sigma'|\eta')}\right\}\COMMA
 \end{align}
where with a slight abuse of notation we denote the probability density of measures with the same letter. Note that in the sequel $\eta=\COP\sigma, \eta'=\COP\sigma'$, if not
otherwise stated. 
A trial state $ \sigma'$ is  generated by the  first level and the simple reconstruction with the transition probability
\begin{equation*}
Q(\sigma,\sigma') = \alpha_{CG}(\eta,\eta')\BARIT{\rho}(\eta,\eta')\mu_r(\sigma'|\eta'),\; \mbox{for $\sigma\neq \sigma' $} \PERIOD
\end{equation*}
 It is easy to check that $Q(\sigma,\sigma')$ satisfies the detailed balance condition, see Definition~\ref{defDB},  with   $ \mu_0(d\sigma) =\BARMZ(d\eta) \mu_r(d\sigma|\eta)$.
i.e. $  Q(\sigma,\sigma') \mu_0(\sigma) = Q(\sigma',\sigma)\mu_0(\sigma') $.
This property lead to the simplified formulation of the second level acceptance probability \VIZ{gfineAcc} from  the one suggested by the Metropolis method  $\alpha_f(\sigma,\sigma')=\min\left\{1,\frac{\mu(\sigma')Q(\sigma',\sigma)}{\mu(\sigma) Q(\sigma,\sigma')}\right\} $.
The probability of moving from a state $\sigma$ to a next $\sigma'$ is $\alpha_f(\sigma,\sigma')\alpha_{CG}(\eta,\eta')\mu_r(\sigma'|\eta')\BARIT{\rho}(\eta,\eta')$. Therefore the method generates a   Markov chain $\{Y_n\}$, starting from an arbitrary initial state $ \sigma_0$,  with transition kernel
\begin{eqnarray}\label{KCG}
&\Kk_{CG}(\sigma,\sigma')=\alpha_f(\sigma,\sigma')\alpha_{CG}(\eta,\eta')\mu_r(\sigma'|\eta')\BARIT{\rho}(\eta,\eta') \quad \text{for } \sigma\neq \sigma' \COMMA  \\
&\Kk_{CG}(\sigma,\sigma)= 1-\int_{\Sigma}\alpha_f(\sigma,\sigma')\alpha_{CG}(\eta,\eta')\mu_r(\sigma'|\eta')\BARIT{\rho}(\eta,\eta')  d\sigma' \nonumber
\PERIOD
\end{eqnarray}

We denote $E=\{\sigma\in\Sigma:\; \mu(\sigma)>0 \} $ and $\tilde E=\{ \sigma \in \Sigma :\;  \mu_0(\sigma)>0 \} $ the support of the microscopic and the {\it proposal}  distributions respectively.

With the following Theorem,  the proof  of which is given in \ref{Math},
we prove that transition kernel $ \Kk_{CG}(\sigma,\sigma')$ satisfies the detailed balance condition, that   ensures the method generates samples from the target measure. Furthermore,  irreducibility and aperiodicity properties are  satisfied that guarantee ergodicity of  $\{ Y_n\} $,  i.e. $\frac{1}{n} \sum_{j=1}^n f(Y_j) $ is a convergent approximation of the averages $ \int f(\sigma) \mu(d\sigma) $ for any $f\in L^1(\mu)$.
\begin{thm}\label{ReversibilityThm}
For every conditional distribution $\bar{\rho}(\eta,\eta') $ on $\BARIT{\Sigma} $,
and $\mu_r(\cdot|\eta)$ on $ \Cspace $,
 \begin{enumerate}[i)]
	\item The transition kernel $\Kk_{CG}(\sigma,\sigma')$ satisfies the detailed balance (DB)
condition with $ \mu(\sigma) $.
	 \item    $\mu(\sigma)$ is a stationary distribution of the chain.
\item If $\bar\rho(\eta,\eta')>0$,  $  \mu_r(\sigma|\eta)>0$ for all $\sigma,\sigma'\in E$ and $E \subset \tilde{E}$ holds, then $\{Y_n\}$ is $\mu$--irreducible.
	 \item  $\{Y_n\}$ is aperiodic.
 \end{enumerate}
\end{thm}

\section{Computational Complexity of Multilevel CGMC}\label{complexity}
The effectiveness of the multilevel CGMC method is a result of the synthesis of the following two arguments. Firstly,  the { \it
computational cost } of a conventional MH method is reduced  by a two--level CGMC method.  However this is not enough to prove that the method can indeed accelerate conventional methods, and a mathematical spectral analysis emerges as a key algorithmic need. Therefore, secondly, Theorem~\ref{main} provides the relation of the two methods equilibration times using spectral arguments.

\subsection{Computational complexity}
An abstract comparison of the computational complexity of   a conventional MH and the two--level CGMC is summarized in  Table \ref{complex}, for sampling a Gibbs measure with Hamiltonian \VIZ{HamSplit}. By computational complexity here we mean the cost of calculating energy differences involved at the acceptance probabilities. The following  analysis is based on the approximate CG strategy with potential splitting   as described in subsection \ref{Decompositions}.  
Let us consider the canonical ensemble where energy difference  is
\begin{eqnarray*}
 \Delta H_N(\sigma, \sigma^x) = (2\sigma(x) -1) \sum_{ y\in \LATT,y\neq x} \left[K(x-y) +    J(x-y)\sigma(y)\right] \PERIOD
\end{eqnarray*}
For potential $J(x-y)$ with interaction range $L$  each particle interacts with a number of $(2L+1)^d $ neighbours and similarly for $K(x-y)$ with range $S$.  Therefore the number of operations necessary   are   $ (2L+1)^d + (2S+1)^d $. Similarly the energy differences appearing in the two--level method are
 \begin{eqnarray*}
 &\Delta \BARHZ_l(\eta, \eta^k) & = \sum_{l\in\LATTC,  k\neq l}\BARIT{J}(k,l) \eta(l)+ \BARIT{J}(0,0)(\eta(k)-1)\COMMA\\
&\Delta H_s(\sigma, \sigma^x) &=  (2\sigma(x) -1)  \sum_{y\in\LATT,   y\neq x} K(x-y)  \sigma(y)\COMMA
\end{eqnarray*}
at the first  and    second level respectively. The compressed interaction potential  $\BARIT{J}(k,l)$ has   the reduced  range  $ L/q $   and the number of operations for calculating    $ \Delta \BARHZ_l(\eta, \eta^k)  $ is $ (2L+1)^d /Q $, while  for $ \Delta H_s(\sigma, \sigma^x) $    is $ (2S+1)^d /Q $.  
The method,  in addition to the reduction of operations due to  range suppression,  exhibits  a computational  reduction  as a result of the   fact  that rejected trials   at the first level will  not be tested at the second, and the calculation of differences $ \Delta H_s(\sigma, \sigma^x) $   is avoided.   A summary of this discussion appears in Table~\ref{complex}.
\begin{table}
\caption{Operations count for evaluating energy differences for $n$
iterations.  The total number of accepted coarse trials $m$ is  the  number of   the second level iterations tested.  }
  \begin{tabular}{lc} 
\hline\noalign{\smallskip}
 		Method 		&  Operations \\
\hline\noalign{\smallskip}
Metropolis Hastings	&   $ n \times (2L+1)^d + n\times  (2S+1)^d $\\ 
Two-level CGMC 	&  $ n \times (2L+1)^d/Q + m  \times  (2S+1)^d$	    \\
\noalign{\smallskip}\hline\noalign{\smallskip}
 \end{tabular}
 \label{complex}
 \centering
 \end{table}
  When the correction terms $ \Delta H_l(\sigma,\sigma^x) - \Delta \BARHZ_l(\eta,\eta^k) $ are present  the additional computational  cost in the second level is small,  dependent on the decay of rate of  $J(x-y)-\BARIT{J}(k,l),\ x,y \in \Sigma_N, x\in C_k, y\in C_l $. 
Since this term is in principle fast decaying,  in implementations we can  neglect interactions with  distance larger than a cut-off range $ L_c< L$ with a small error, that will contribute   further in the  computational time reduction. 
 
\subsection{Comparison of equilibration times}
A direct estimation of rate of convergence of a Markov chain generated by a Monte Carlo method is model dependent and in general intractable. Thus for the purpose of this work it is natural to
study  this property as a comparison of the proposed method with a
conventional method. Theorem~\ref{main} provides such a comparison of the spectral gap between the conventional MH method and the two-level CGMC method, for sampling a measure $ \mu(d\sigma)$ on a countable state space $\Sigma$. This comparison is summarised in inequality \VIZ{gap}  proving that their relation, in terms  the spectral gap,  is controlled by the approximation $ \BARMZ(d\eta)$ of marginals $ \BARIT{\mu}(d\eta) = \mu \circ \COP^{-1}(d\eta)$.

For a discrete time Markov chain $\{X_n\}$ with
transition kernel $\Kk$ and stationary distribution $\mu$, the mixing time $
\tau $ is defined as
\begin{equation*}
\tau:= min_n \left\{ \forall \sigma\in \Sigma : \| \Kk^n(\sigma,\cdot) -
\mu(\cdot)\|_{TV} \le \frac14\right\} \PERIOD
   \end{equation*}
For two probability measures $ \mu, \nu $ the total variation norm is $\| \mu - \nu\|_{TV} =\frac{1}{2}\sum_{\sigma}| \mu(\sigma) - \nu(\sigma)|$.
Bounds of the total variation norm appearing can be given  in terms of $\Kk$'s
spectral gap $\lambda$, for example for a  reversible kernel holds \cite{DSC}
\begin{equation}\label{bound}
2||\Kk^n(\sigma,\cdot) -\mu ||_{TV} \le \frac{1}{ \mu(\sigma)^{1/2}}(1-\lambda)^n \PERIOD
\end{equation}  
The spectral gap of a kernel $\Kk$ is defined by $\lambda(\Kk)=\min
\left\{\frac{\mathcal{E}(f,f)}{\VAR(f)}; \VAR(f)\neq 0\right\}$ with the
Dirichlet form $\mathcal{E}(f,f)=\frac12
\sum_{\sigma,\sigma'}|f(\sigma)-f\sigma')|^2\Kk(\sigma,\sigma')\mu(\sigma)$ and
the variance $ \VAR(f)=\frac12
\sum_{\sigma,\sigma'}|f(\sigma)-f(\sigma')|^2\mu(\sigma')\mu(\sigma)$.
 In view of \VIZ{bound}, one can say that between two algorithms producing Markov chains with identical equilibrium distributions {\it better} in terms of speed of convergence is the one with the {\it larger spectral gap}.

Let $\mathcal{E}(\Kk_{CG}), \mathcal{E}(\Kk_{c})$ denote the Dirichlet forms
and $\lambda(\Kk_{CG}), \lambda(\Kk_{c})$ the spectral gap corresponding to the
two-level CGMC and the classical Metropolis  transition kernels $\Kk_{CG}(\sigma,\sigma')$ \VIZ{KCG} and $\Kk_c(\sigma,\sigma') $ \VIZ{KC} respectively.


\begin{thm}\label{main}
Let $\rho(\sigma,\sigma') $ be a symmetric proposal transition probability for
the conventional MH algorithm and $\bar{\rho}(\eta,\eta') $ a symmetric
proposal transition probability on the coarse space $ \bar{\Sigma}$ for the
two-level CGMC algorithm, then for any  conditional probability $\mu_r(\sigma|\eta)$
 \begin{equation}\label{gap}  
\mathcal{A}\underline{\gamma}\lambda(\Kk_c) \le \lambda(\Kk_{CG})\le
\bar{\gamma}\lambda(\Kk_c)
\end{equation}
where $\mathcal{A}=\inf_{\sigma,\sigma'}\{\mathcal{ A}(\sigma,\sigma')\}$
and $\underline{\gamma}>0, \bar{\gamma}>0 $ such that  
$ \underline{\gamma}\le \Bb(\sigma,\sigma')\le \bar{\gamma}$, with $\mathcal{
A}(\sigma,\sigma')$ and $\Bb(\sigma,\sigma')$ defined in lemma \VIZ{KClemma}.
  \end{thm}
 \noindent
{\bf Remarks.} 
\begin{enumerate}
\item Existence of finite and positive values of  $\mathcal{A}=\inf_{\sigma,\sigma'}\{\mathcal{
A}(\sigma,\sigma')\}  $ is ensured for the models studied in this work.  Consider the Gibbs measure  $  \mu_{N,\beta}(d\sigma) \propto \exp\{-\beta H_N(\sigma)\}$ and $ 
\BARMZ_{M,\beta}(d\eta) \propto \exp\{-\beta \BARHZ (\eta)\}$ as defined
in Section \ref{SLS}.  Let  $ J(r), s.t. \sum_r  J(r) = J^*<\infty $ with the compressed  interactiosn $\bar J $ as defined in \VIZ{Hzero}. All possible values of $\mathcal{A}(\sigma,\sigma') $, \VIZ{DefA}, are 
$\mathcal{A}(\sigma,\sigma')= 1$, 
\begin{equation*}
 \mathcal{A}(\sigma,\sigma')= \min\{ \EXP{-\beta\Delta \BARHZ(\eta,\eta')} ,\EXP{-\beta\Delta \BARHZ (\eta',\eta)}\} \ge e^{-\beta \bar J^*} \PERIOD
\end{equation*}
and
\begin{eqnarray*}
\mathcal{A}(\sigma,\sigma')&=&  
\min\{ e^{-\beta(\Delta H_N(\sigma,\sigma') -\Delta \BARHZ(\eta,\eta'))
},e^{-\beta(\Delta H_N(\sigma',\sigma) -\Delta \BARHZ(\eta',\eta))}\} \\
&=& 1+ \BIGO\left(\frac{Q}{(2L+1)^d}\right) 
\PERIOD
 \end{eqnarray*}
This is a result of the known  estimate  Lemma (2.3) \cite{KPRT},
\begin{equation*}
 \Delta H_N(\sigma,\sigma') -\Delta \bar
H^{(0)}(\eta) = \BIGO\left(\frac{Q}{(2L+1)^d}\right)
\end{equation*}
This  proves also that values of $\Aa(\sigma,\sigma')$ are close to $1$ controlled by the approximation parameter $ Q$.
 \item Term $ \Bb(\sigma,\sigma')$ depends only  on the difference
between the proposal kernels in the two methods. For example if  $ \mu_r(\sigma|\eta)$ is  chosen such that 
$ \rho(\sigma,\sigma') = \bar\rho(\eta,\eta')\mu_r(\sigma'|\eta')$ for
all $ \sigma, \sigma'\in \Sigma$
then $ \Bb(\sigma,\sigma') \equiv 1$, and inequality \VIZ{gap}
depends only on $\Aa$, that is
  \begin{equation*}  
\mathcal{A} \lambda(\Kk_c) \le \lambda(\Kk_{CG})\le
\lambda(\Kk_c)\PERIOD
\end{equation*}
\end{enumerate}

Numerically the statement of the Theorem is revealed by a comparison of the
average acceptance probabilities in examples, see Figures \ref{fig:2.1}. We
expect that best approach in terms of equilibration rate is a rejection free
method for sampling the CG distribution, which is an approach that we elaborate
with in a forthcoming work.

To prove Theorem~\ref{main} we need the following lemmata.
\begin{lem}\label{KClemma}
For symmetric proposal transition kernels 
$\rho(\sigma,\sigma')=\rho(\sigma',\sigma)$ and
$\bar{\rho}(\eta,\eta')=\bar{\rho}(\eta',\eta)$,  
for all $\sigma,\sigma'\in \Sigma$ with $\sigma\neq\sigma'$, 
\begin{equation*} 
\Kk_{CG}(\sigma,\sigma')=\mathcal{A}(\sigma,\sigma')\mathcal{B}(\sigma,\sigma')
\Kk_{c}(\sigma,\sigma')
  \end{equation*}
  where
\begin{equation}\label{DefB}
\mathcal{B}(\sigma,\sigma') =\frac{\bar{\rho}(\eta,\eta')
}{\rho(\sigma,\sigma')} \begin{cases}
\mu_r(\sigma'|\eta'), \; &\mbox{ if  $\alpha_f(\sigma,\sigma')=1$}\\
\mu_r(\sigma|\eta), \; &\mbox{ if  $ \alpha_f(\sigma,\sigma')<1$}
\end{cases}  \PERIOD
\end{equation}
Furthermore we define the subsets
\small{
\begin{eqnarray}\label{sets}
&C_1=\!\left\lbrace (\sigma,\sigma')\!\in \!\Sigma\!\times\!\Sigma\!: \left\lbrace
\!\alpha<1, \alpha_{CG}<1, \alpha_f<1\! \right\rbrace  \text{ or } 
\left\lbrace \!\alpha=1, \alpha_{CG}=1, \alpha_f=1\! \right\rbrace
\right\rbrace\\
&C_2=\!\left\lbrace (\sigma,\sigma')\!\in \! \Sigma\!\times\!\Sigma\!:\left\lbrace
\!\alpha=1, \alpha_{CG}<1, \alpha_f=1 \! \right\rbrace\text{ or }
\left\lbrace \! \alpha<1, \alpha_{CG}=1, \alpha_f<1 \!\right\rbrace
\right\rbrace\nonumber\\
&C_3=\!\left\lbrace (\sigma,\sigma')\!\in \! \Sigma\!\times\!\Sigma\!: \left\lbrace
\!\alpha=1, \alpha_{CG}=1, \alpha_f<1 \!\right\rbrace  \text{ or } 
\left\lbrace \!\alpha<1, \alpha_{CG}<1, \alpha_f=1 \!\right\rbrace
\right\rbrace\nonumber\\
&C_4=\!\left\lbrace (\sigma,\sigma')\!\in \!\Sigma\!\times\!\Sigma\!:
\left\lbrace \!\alpha<1, \alpha_{CG}=1, \alpha_f=1 \!\right\rbrace \text{
or } \left\lbrace \! \alpha=1, \alpha_{CG}<1, \alpha_f<1\! \right\rbrace
\right\rbrace \nonumber
\end{eqnarray}
}
and
\begin{equation}\label{DefA}
\mathcal{A}(\sigma,\sigma')=
\begin{cases}
	    1,\; \;  &\mbox{ if $ (\sigma,\sigma')\in  C_1$} \\
\min\{
\frac{\BARMZ(\eta')}{\BARMZ(\eta)},\frac{\BARMZ(\eta)}{\BARMZ(\eta')}\}, \; \;       & \mbox{ if  $(\sigma,\sigma')\in  C_2$}  \\
\min\{\frac{\mu(\sigma')\BARMZ(\eta)}{\mu(\sigma)\BARMZ(\eta')},
\frac{\mu(\sigma)\BARMZ(\eta')}{\mu(\sigma')\BARMZ(\eta)}\}, 
\; \;   &\mbox{ if  $(\sigma,\sigma')\in  C_3$}   \\
\min\{ \frac{\mu(\sigma')}{\mu(\sigma)},\frac{\mu(\sigma)}{\mu(\sigma')}\}, 
\; \;   &\mbox{ if  $(\sigma,\sigma')\in  C_4$} 
\end{cases}
\end{equation} 
\end{lem}

\noindent
\pf 
Recall that
$\Kk_{c}(\sigma,\sigma')=\alpha(\sigma,\sigma')\rho(\sigma,\sigma')$, for $\sigma\neq\sigma' $, which for the conventional MH method and the hypothesis of symmetry of $\rho$
becomes $\Kk_{c}(\sigma,\sigma')=\min\left\{1,\frac{\mu(\sigma')}{\mu(\sigma)}
\right\}\rho(\sigma,\sigma')$. Similarly the two--level CGMC   kernel becomes
\begin{eqnarray*} 
\Kk_{CG}(\sigma,\sigma')=&&\min\left\{1,\frac{\mu(\sigma')\BARMZ(\eta)\mu_r(\sigma|\eta)}{\mu(\sigma)\BARMZ(\eta') \mu_r(\sigma'|\eta')}\right\}\\
&&\times\min\left\{1,\frac{\BARMZ(\eta')}{\BARMZ(\eta)}\right\}\bar{\rho}(\eta,\eta')\mu_r(\sigma'|\eta')\PERIOD
  \end{eqnarray*}  
The proof is  based on a case study on the values of the acceptance
probabilities $ \alpha(\sigma,\sigma') = \min\left\{1,\frac{\mu(\sigma')}{\mu(\sigma)} \right\}$,
$\alpha_f(\sigma, \sigma') = \min\left\{1,\frac{\mu(\sigma')\BARMZ(\eta)\mu_r(\sigma|\eta)}{\mu(\sigma)\BARMZ(\eta')\mu_r(\sigma'|\eta')}\right\}$
and $\alpha_{CG}(\eta, \eta')=
\min\left\{1,\frac{\BARMZ(\eta')}{\BARMZ(\eta)}\right\} $.
In the following we also use extensively the symmetry property of $\rho$ and
$\bar\rho $.  All possible combinations on their values are categorized in  the sets
$C_i, i=1,\dots,4$ defined in \VIZ{sets}.
For $(\sigma, \sigma') \in C_1$ and $\alpha(\sigma, \sigma')<1$,
$\alpha_f(\sigma, \sigma')<1$ and $\alpha_{CG}(\eta, \eta')<1$,
$ \Kk_{c}(\sigma,\sigma')= \frac{\mu(\sigma')}{\mu(\sigma)}
\rho(\sigma,\sigma') $ and
\begin{equation*}
\Kk_{CG}(\sigma,\sigma')=
\frac{\mu(\sigma')\BARMZ(\eta)}{\mu(\sigma)\BARMZ(\eta')}
\frac{\BARMZ(\eta')}{\BARMZ(\eta)}\bar{\rho}(\eta,\eta')\mu_r(\sigma'|\eta')=
\frac{\mu(\sigma')}{\mu(\sigma)}\bar{\rho}(\eta,\eta')\mu_r(\sigma'|\eta')\COMMA
\end{equation*}
and we derive their relation  $\Kk_{CG}(\sigma,\sigma')=\frac{\bar{\rho}(\eta,\eta')\mu_r(\sigma|\eta)}
{\rho(\sigma,\sigma')}\Kk_c(\sigma,\sigma')$.
With similar simple calculations, that we omit here, for all sub--cases that
are encountered in subsets $C_i, i=1,\dots 4 $ we have, for $(\sigma, \sigma') \in C_1$ 
\begin{equation*}
\Kk_{CG}(\sigma,\sigma')=
\Bb(\sigma,\sigma')
\Kk_{c}(\sigma,\sigma')\PERIOD
\end{equation*}
with 
\begin{equation*}
\mathcal{B}(\sigma,\sigma') =\frac{\bar{\rho}(\eta,\eta')
}{\rho(\sigma,\sigma')} \begin{cases}
\mu_r(\sigma'|\eta'), \; &\mbox{ if  $\alpha_f(\sigma,\sigma')=1$}\\
\mu_r(\sigma|\eta), \; &\mbox{ if  $ \alpha_f(\sigma,\sigma')<1$}
\end{cases} \PERIOD
\end{equation*}
For $(\sigma, \sigma') \in C_2$  such that 
$\alpha(\sigma, \sigma')=1$, $\alpha_{CG}(\eta, \eta')<1,
\alpha_f(\sigma, \sigma')=1$,
\begin{equation*}
\Kk_{CG}(\sigma,\sigma')= \frac{\BARMZ(\eta')}{\BARMZ(\eta)}
\frac{\bar{\rho}(\eta,\eta')\mu_r(\sigma'|\eta')}{\rho(\sigma,\sigma')}
\Kk_{c}(\sigma,\sigma')\COMMA
\end{equation*}
 and for
$\alpha(\sigma, \sigma')<1$, $\alpha_{CG}(\eta, \eta')=1$,
$\alpha_f(\sigma, \sigma')<1$,
\begin{equation*}
\Kk_{CG}(\sigma,\sigma') =
\frac{\BARMZ(\eta)}{\BARMZ(\eta')}
\frac{\bar{\rho}(\eta,\eta')\mu_r(\sigma|\eta)}{\rho(\sigma,\sigma')}
K_{c}(\sigma,\sigma')\COMMA
\end{equation*}
 that we can summarize to 
\begin{equation*}
\Kk_{CG}(\sigma,\sigma') =
\min\left\{\frac{\BARMZ(\eta)}{\BARMZ(\eta')}, \frac{\BARMZ(\eta')}{\BARMZ(\eta)}\right\} \Bb(\sigma,\sigma')K_{c}(\sigma,\sigma')\COMMA
\end{equation*}
 since for example in the first case $\frac{\BARMZ(\eta')}{\BARMZ(\eta)}< 1 < \frac{\BARMZ(\eta)}{\BARMZ(\eta')}$ and inversely for the second.
Following the same reasoning for $(\sigma, \sigma') \in C_3$,
\begin{equation*}
\Kk_{CG}(\sigma,\sigma')=
\min\left\{  \frac{\mu(\sigma')\BARMZ(\eta)}{\mu(\sigma)\BARMZ(\eta')}, \frac{\mu(\sigma)\BARMZ(\eta')}{\mu(\sigma')\BARMZ(\eta)} \right\}
\Bb(\sigma,\sigma') \Kk_{c}(\sigma,\sigma')\COMMA
\end{equation*}
and for  $(\sigma, \sigma') \in C_4$
\begin{equation*}
\Kk_{CG}(\sigma,\sigma')= \min\left\{ \frac{\mu(\sigma') }{\mu(\sigma) }, \frac{\mu(\sigma) }{\mu(\sigma') } \right\}
\Bb(\sigma,\sigma')\Kk_{c}(\sigma,\sigma')\PERIOD
\end{equation*}
All these steps prove, the following relation of transition kernels generated by  Algorithms~\ref{MH} and \ref{mCGMC},
\begin{equation*}
\Kk_{CG}(\sigma,\sigma')= \Aa(\sigma,\sigma') \mathcal{
B}(\sigma,\sigma') \Kk_{c}(\sigma,\sigma')\COMMA
\end{equation*}  
with $ \mathcal{ A}(\sigma,\sigma') $ and $ \Bb(\sigma,\sigma') $ defined in
\VIZ{DefA} and \VIZ{DefB}.
\hfill $ \square$

The proof of Theorem~\ref{main} is based on the application of Lemma~
\ref{KClemma} and Lemma 3.3 in the work of P. Diaconis and L.
Saloff-Coste \cite{DSC} that is stated here for completeness.
\begin{lem}\label{spectral}
Let $\Kk, \mu $ and $\Kk', \mu' $ be two Markov chains on the same finite set
$X$. Assume that there exist $A, a>0$ such that
\begin{equation*}
\Ee'\le A\Ee,\; \;   a\mu\le  \mu' \COMMA
\end{equation*}
then
$$ \lambda'\le \frac{A}{a}\lambda\PERIOD$$
\end{lem}
We conclude with the proof of Theorem~\ref{main}.
 \pf
We compare the Dirichlet forms $\Ee(\Kk_{CG}), \Ee(\Kk_{c})$ using the
definition of Dirichlet form and applying Lemma (\ref{KClemma}). By the definition of $ \Aa(\sigma,\sigma')$ for all $\sigma,\sigma'\in \Sigma$ holds $0<\Aa(\sigma,\sigma')\le 1$ and assume that $\Aa:=\inf_{\sigma,\sigma'}\{\Aa(\sigma,\sigma')\} >0$.
Then    by Lemma \ref{KClemma}  
\begin{equation*} 
\inf_{\sigma,\sigma'}\{\Aa(\sigma,\sigma')\}
\Bb(\sigma,\sigma')\Kk_c(\sigma,\sigma')\le \Kk_{CG}(\sigma,\sigma')\le
\Bb(\sigma,\sigma')\Kk_c(\sigma,\sigma')\PERIOD
 \end{equation*}
Let $\underline{\gamma}>0$ and $\bar{\gamma}>0$ such that 
$\underline{\gamma}\le \Bb(\sigma,\sigma')\le \bar{\gamma}$.
Then 
\[ \inf_{\sigma,\sigma'}\{\Aa(\sigma,\sigma')\}\underline{\gamma}
\Kk_c(\sigma,\sigma')\le\Kk_{CG}(\sigma,\sigma')\le
\bar{\gamma}\Kk_c(\sigma,\sigma')\PERIOD\]
Recalling the definition of Dirichlet form for kernel $ \Kk_{CG}$, 
\[
\Ee_{CG}(f,f)=\frac12
\sum_{\sigma,\sigma'}|f(\sigma)-f(\sigma')|^2\Kk_{CG}(\sigma,\sigma')\mu(\sigma)\COMMA
\]
and the above relation we can write 
$\inf\{\Aa(\sigma,\sigma')\}\underline{\gamma}
\mathcal{E}_c\le \mathcal{E}_{CG}\le \bar{\gamma}\mathcal{E}_c $.
Application of Lemma \ref{spectral},  for which here $\mu'\equiv \mu $, thus $ a=1$ and 
 $A = \inf_{\sigma,\sigma'}\{\Aa(\sigma,\sigma')\} \underline\gamma$ for the left hand side inequality, and $A =  \bar{\gamma}  $ for the right hand side,   gives the  relation of spectral gaps
\begin{equation*}
\inf_{\sigma,\sigma'}\{\Aa(\sigma,\sigma')\}\underline{\gamma}\lambda(\Kk_c) \le
\lambda(\Kk_{CG})\le \bar{\gamma}\lambda(\Kk_c) \PERIOD
\end{equation*} 
 \hfill $\square $


\section{Benchmark examples for the canonical ensemble}\label{applications}
\subsection{Combined Ising and Curie Weiss model}\label{benchmark}
We consider a benchmark example of competing short- and long- range
interactions, where constants $ \mathcal{A}$, and $\underline{\gamma},
\bar\gamma $ appearing in Theorem~\ref{main} are calculated explicitly.
Furthermore analytical expressions for the free energy in the thermodynamic
limit, $N\to \infty$, are known \cite{KAR}. The energy of the system at
configuration $\sigma \in \Sigma_N$ is
\begin{eqnarray*}
H_N(\sigma) &=&
-\frac{K}{2}\sum_{x\in\LATT}\sum_{   |x-y|=1}\sigma(x)\sigma(y)
-\frac{J}{2N}\sum_{x\in\LATT}\sum_{  y\neq x}
\sigma(x)\sigma(y)-h\sum_{x\in\LATT}\sigma(x)\\
&=&  H_s(\sigma)+ H_l(\sigma) -h\sum_{x\in\LATT}\sigma(x)\PERIOD
\end{eqnarray*}
  The interactions involved in $H_s(\sigma)$ are of the nearest-neighbour type with
strength $K$, while $H_l(\sigma) $ represents a mean--field approximation of a potential $J(x-y)$  as in  \VIZ{microJL} for example, averaged over all lattice sites. The coarse grained
Hamiltonian $\BARH_l(\eta) $ is exact, equal to the microscopic energy
$H_l(\sigma) $. Indeed,  for any coarse graining level $ Q$ 
 \begin{equation*}
H_l(\sigma)=\BARH_l(\eta)=-\frac{J}{2N}\sum_{k\in\LATTC}\sum_{l\in\LATTC 
l\neq k}\eta(k)\eta(l) -\frac{J}{2N}\sum_{k\in\LATTC}
\eta(k)\left(\eta(k)-1\right)\COMMA
\end{equation*}
where  $ \eta(k) =\sum_{x\in C_k} \sigma(x), k\in \LATTC $.
No coarse graining approximation error is involved from compressing long range
interactions, that allows us to study the effect of  the splitting of short-- and long-- range interactions.  We consider the canonical ensemble with the spin-flip dynamics. 

 The detailed calculations on the application of Theorem~\ref{main} are given  in   \ref{appendix1}, where we prove that the constants appearing in inequality \VIZ{gap} are 
\begin{eqnarray*}
&\underline{\gamma} = \bar{\gamma} =1,\text{ and } \\
&\Aa =
\begin{cases}
\min\{e^{-|J|},e^{-2|K|}\}\; &\mbox{for $K\neq 0$,} \\
1                                     \;  &\mbox{for  $K=0$.} \\
\end{cases}
\end{eqnarray*}
The relation of the spectral gap of the methods is,  for $K\neq 0$,
\begin{equation*} 
\min\{e^{-|J|},e^{-2|K|} \} \lambda(\Kk_{c})\le \lambda(\Kk_{CG})\le
\lambda(\Kk_{c})\PERIOD
\end{equation*}
For $ K=0$   Theorem~\ref{main}   verifies  the equivalence of the two methods $\lambda(\Kk_{CG})=\lambda(\Kk_{c}) $ as is expected.
The dependence on the splitting of energy is revealed, showing that the   governing parameter is the  relative strength of the short and long range interactions $J$ and $K$.

\begin{figure}[ht]
\centering
\includegraphics[angle=-0,width=0.45\textwidth,
height=0.3\textheight]{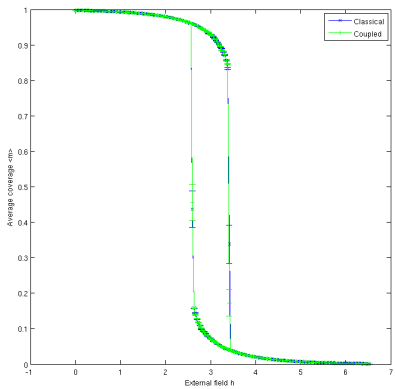}
\caption{ Comparison of hysteresis of the microscopic  
with the two--level process for model \ref{benchmark}.  $ N= 64\times64,\; q=4,\; K= 1,\; J= 5$.
}\label{fig:1}
\end{figure}

Here we  present simulation results for two dimensional lattice systems,  while in \cite{KKP} one can find a detailed numerical study for  one--dimensional  systems.
Table \ref{costL} shows a comparison of the computational time between the MH  and the two--level CGMC method for a   lattice $N=16 \!\times \!16$ and coarsening parameter $Q=q\!\times\! q$, for $q=4$ and $ q=8$. The results  depict a reduce of computational effort at a rate close to $ \BIGO(Q)$. This rate  is not exactly $ \BIGO(Q) $ because of the additional computational time necessary for  implementing the  local uniform reconstruction. The hysteresis diagram is indicated in Figure \ref{fig:1},  for the average total coverage $\AVG{c}$ first upon increasing the field $h$ from low values and then decreasing it from high values.   The total coverage is $c(\sigma)=N^{-1}\sum_x\sigma(x)$ and average  quantities $\AVG{c}$ are  computed after equilibration in the Monte Carlo sampling, using $1000$ samples.
The tests demonstrate that the two--level CGMC method predicts correctly the phase transition regime, compared to the conventional MH, with a reduced computational cost.
\begin{table}[ht]
\caption{CPU cost comparisons for different resolutions in hysteresis simulations. $ 
K = 1,\; J= 5,\; h\in [0,6],\; N = 16\!\times\!16$.}
\centering  
\begin{tabular}{  l c   }  
\hline\noalign{\smallskip} 
Method &   CPU (min)        \\
[0.5ex] \noalign{\smallskip}\hline\noalign{\smallskip}
MH q=1        &  327             \\
Two--level CGMC q=4 &  28          \\
Two--level CGMC q=8 &   9   \\
\hline  
\end{tabular}
\label{costL}
\end{table}

%
%
\subsection{Kac type interactions}\label{nontrivial}
\begin{figure}[htp]
 \centering	 
\includegraphics[angle=-0,width=0.45\textwidth,
height=0.3\textheight]{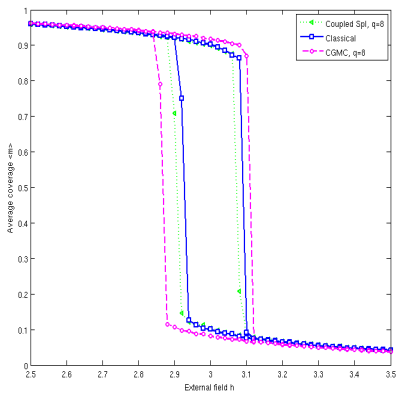}
  \caption{ Comparison of hysteresis of the microscopic  
with the compressed (CGMC) and the  two--level (ML CGMC) process for model \ref{nontrivial}. $\beta J_0 =6,$ $N= 512$, $q= 8$, $S=1$.  }
\label{fig:2}        
\end{figure}

\begin{figure}[htp]
 \centering	 
\subfloat[ ]{{\label{fig:2.1a}}{\includegraphics[scale=0.4
]{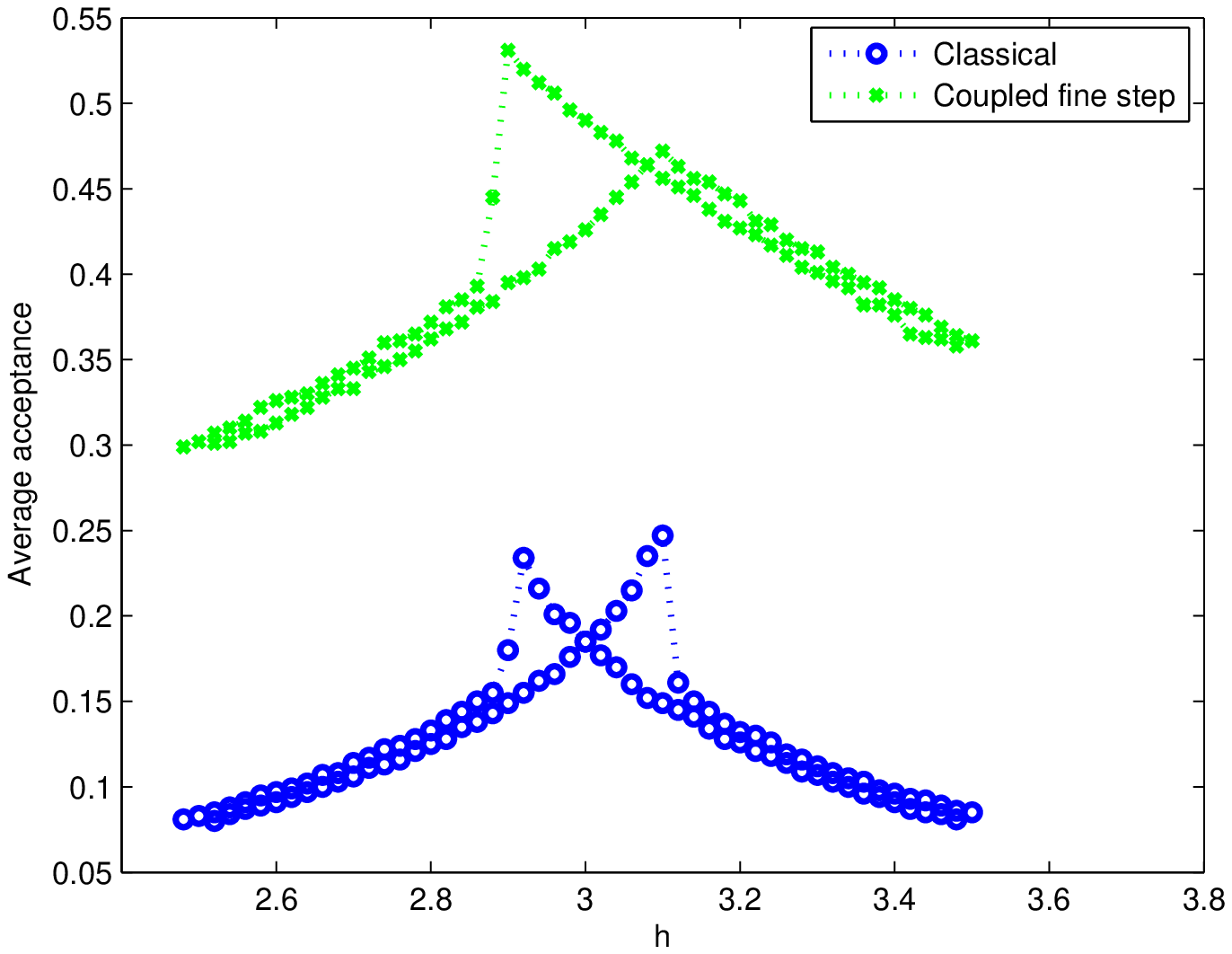}} {\unitlength=1bp}}
\subfloat[ ]{{\label{fig:2.1b}}{\includegraphics[scale=0.4
]{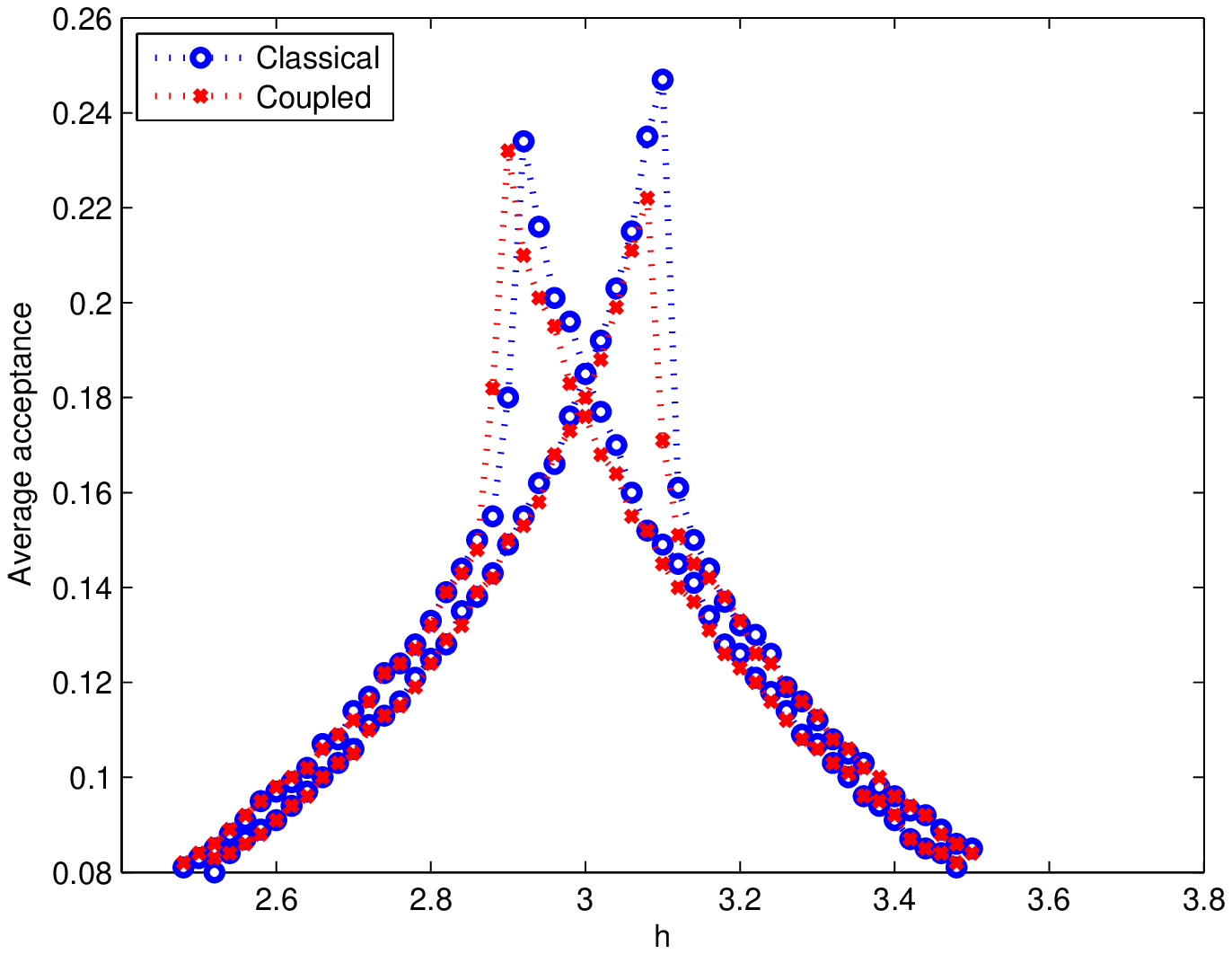}} {\unitlength=1bp}}
\caption{ (a)   The second step ML CGMC and the conventional MH average acceptance probabilities in   hysteresis simulations, Figure \ref{fig:2} showing increase of   acceptances in  a Metropolis - Hastings method.
(b) Comparison of  the total  ML CGMC and the conventional MH average acceptance showing equivalent equilibration rates.}
\label{fig:2.1}        
\end{figure}

 Known error estimates indicate that the potential decay is one of the parameters controlling the approximation error in coarse graining techniques. The  effect of potential singularities on the coarsening error and   improving strategies    encountering  multi-body interactions and/or  potential splitting  has been studied in \cite{AKPR}.
The  example  tested in this section is chosen to study the  potential  splitting strategy with the  proposed method,  that  is shown to improve  direct CGMC results.
We consider  a long range Kac type  interaction potential $ J(r)= N^{-1}V(r/N),\; r=|x-y|>0, x,y\in \LATT$ with an algebraic decay
\begin{equation*}
V(r)=\frac{v}{r^{3/2}},\,  \text{ if } r > 0 \COMMA
\end{equation*} 
where constant $v $ is chosen  to ensure the conservation of the  total mass $J_0=\int J(r)dr $.  The potential splitting strategy  is applied by  decomposing the  interacting potential into a short-range   $ J_s(r)$, as well as a long-range  $ J_l(r)$, defined  by
 \begin{equation*}
J_s(r)=\begin{cases}
J(r)\; &\mbox{for $0< r \le S$,} \\
0  \;  &\mbox{for  $S < r$}, \\
\end{cases}
\end{equation*}
and  
\begin{equation*}
J_l(r)= J(r) - J_s(r) \PERIOD
\end{equation*}
This splitting defines the energy of teh system in the form $ H_N(\sigma) = H_s(\sigma) + H_l(\sigma)$ as in \VIZ{HamSplit} with $ K(x-y) = J_s(|x-y|)$ and $ J(x-y) = J_l(|x-y|)$. In Figure \ref{fig:2} and Table \ref{costJP} we present simulation results  in the canonical ensemble derived from the three methods, the conventional MH with interaction potential $ J(|x-y|), x,y\in \LATT$ and,   the two--level CGMC with interactions $\bar J_l(k,l), k,l\in \LATTC  $ and $ J_s(|x-y|)$, sampling  on the microscopic space $\Sigma$, and the CGMC with interaction potential $\bar J$ sampling  on the coarse space $ \BARIT{\Sigma}_M$. It is observed that, for a fixed coarsening parameter $ Q$,   the two--level method reduces the coarse graining error  of  CGMC simulations, since the most singular part of $ J(r)$, $ J_s(r)$,  is treated in the microscopic space and coarsening is applied only to its   fast decaying part $ J_l(r)$, Figure \ref{fig:2}. 
  Figure \ref{fig:2.1} demonstrates the  increase of acceptance rate of Metropolis sampling when viewing the proposed method as an improving  the  propoposed samples. Overall, as is 
stated in Theorem~\ref{main}, the method's acceptance rate is comparable to  the conventional 
method.  Additionally, a  reduction of the computational time compared to the conventional sampling method  is achieved.
The error appearing in Table \ref{costJP} is the $l^2$ distance of the average total
coverage  $\AVG{c}$, of the indicated method, from  the conventional MH result, as a function of the external field $h$.
\begin{table}[ht] 
\caption{ CPU cost and error  comparisons for different resolutions and methods. $N=512\;,\beta J_0=6,\;S=1$.}  
\centering  
\begin{tabular}{l c c}  
\hline\noalign{\smallskip} 
  Method				&  CPU time&   Error  \\
 [0.5ex] \noalign{\smallskip}\hline\noalign{\smallskip}
	MH    q=1			&     252	&    0        \\
	CGMC q=8			&     47	&    0.78        \\ 
	Two--level CGMC q=8&     57	&    0.65        \\
	CGMC q=64		&     5	&    2.10        \\
	Two--level  CGMC q=64 & 6	&    0.81         \\
\hline  
\end{tabular}
\label{costJP}
\end{table}
 
%
%

\section{Nanopattern formation in heteroepitaxy -- Microcanonical
ensemble}\label{patterns}
We study a  {\it
nanopattern formation} problem and show how the proposed method can provide microscopic information,   benefiting from the low computational cost of the coarse graining  technique while refining the error introduced. The  system studied here   is  characterized by the interplay of short ranged  attraction and long ranged repulsion interactions  between particles. 
Such an interplay can lead to the formation of  patterns, such as discs   and  stripes appearing for example  in heteroepitaxy.
  
The energy of the system is given by $ H_N(\sigma)$, as in \VIZ{GenHam}, with isotropic interaction potential
\begin{equation*}
J(r)=J_0\left( e^{-(r/r_a)^2} -\chi e^{-(r/r_r)^2} \right),\; r=|x-y|>0, x,y\in\LATT\PERIOD
\end{equation*}
$J_0$ is the strength of the potential, $r_a$ and $r_r$ are the dimensionless length scales of attraction and repulsion, respectively, and $\chi$ is the repulsion strength parameter.
A  study of the  kinetic phase diagrams and application of CGMC methods for this system is given in \cite{D1}.  
The underlying dynamics describing the surface diffusion of particles considered are spin-exchange, obeying the exclusion principle and conserving the order parameter microcanonically. 
The order parameter here is the total coverage $ c_0= N^{-1}\sum_{x\in \LATT} \sigma(x) $ representing the   number of occupied sites on the lattice.

\begin{figure}[ht]
 \centering
\includegraphics[angle=-0,width=0.45\textwidth,
height=0.3\textheight]{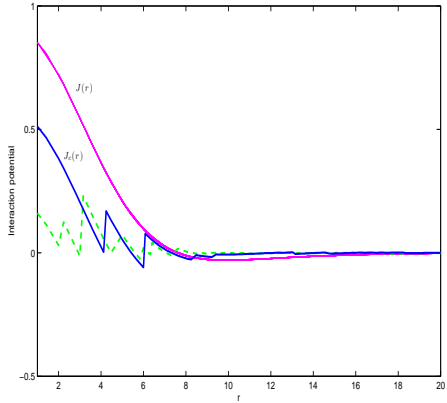} \caption{ Correcting compressed interactions. }
\label{fig:3}        
\end{figure}

\begin{figure}[htp]
         \centering
\subfloat[
MH]{{\label{fig:discs-tr}}{\includegraphics[scale=0.25]{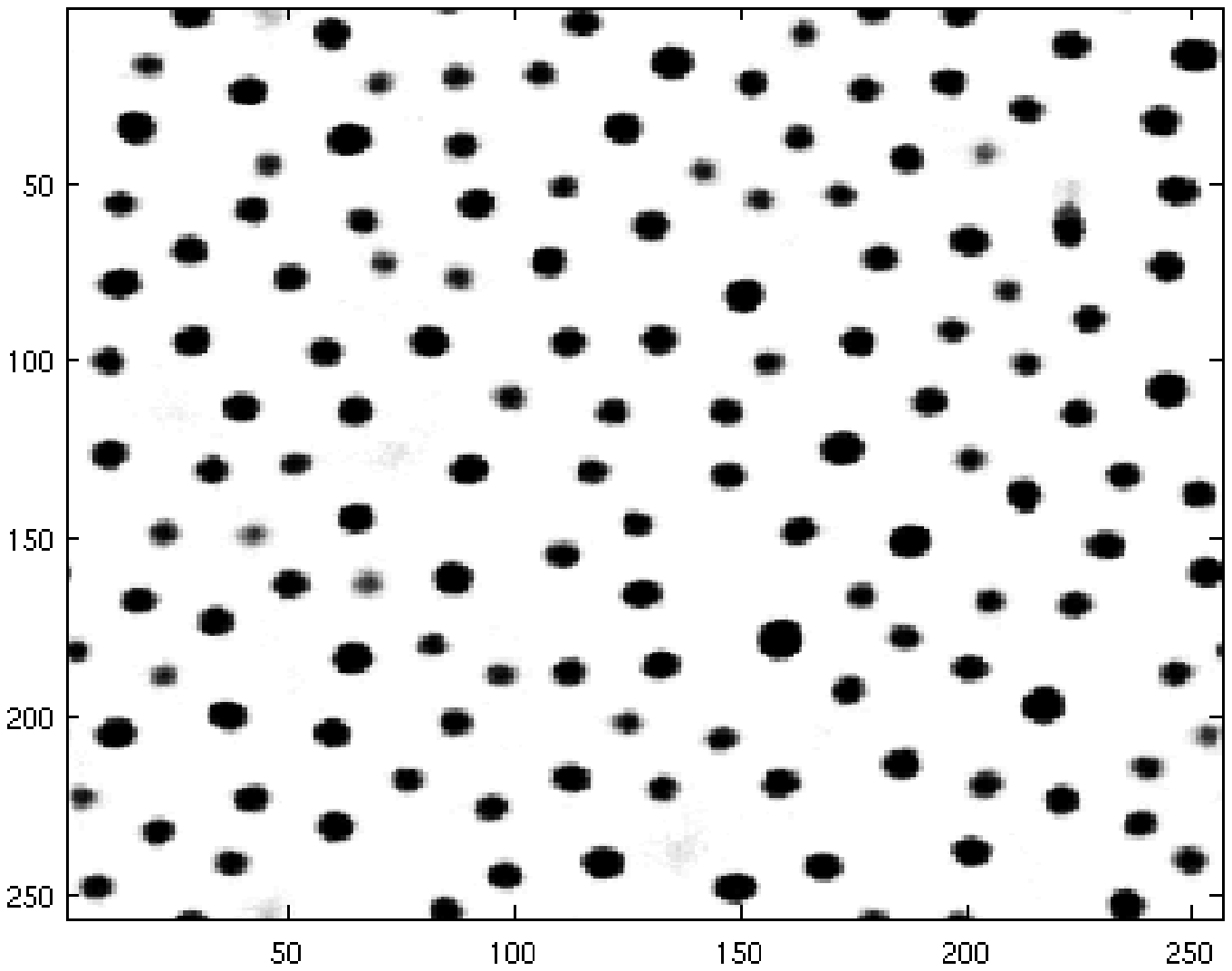}}
{\unitlength=1bp}}
\subfloat[CGMC,
q=8]{{\label{fig:discs-cg}}{\includegraphics[scale=0.25]{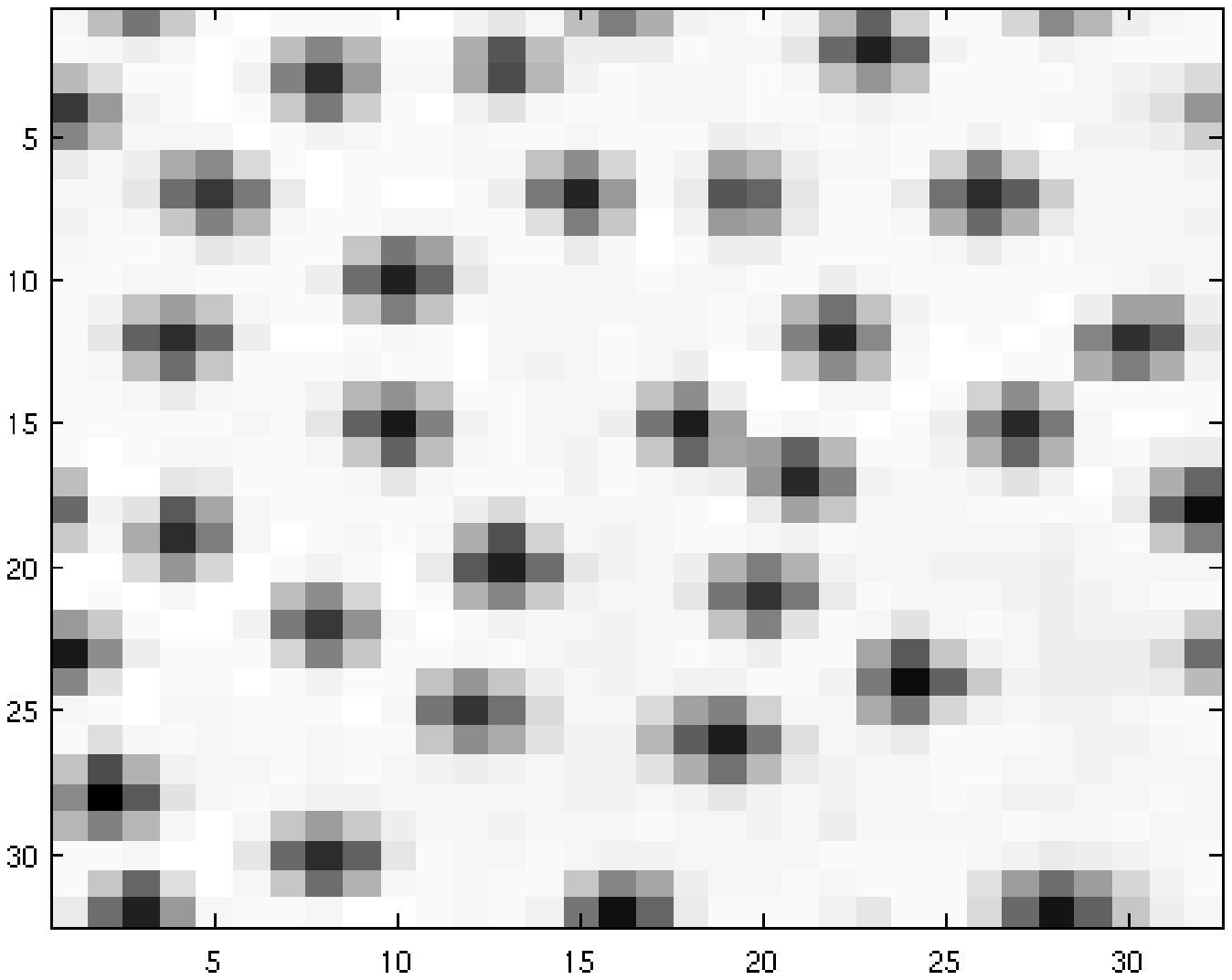}}{\unitlength=1bp}}
\subfloat[Two--level,
q=8]{{\label{fig:discs-co}}{\includegraphics[scale=0.25]{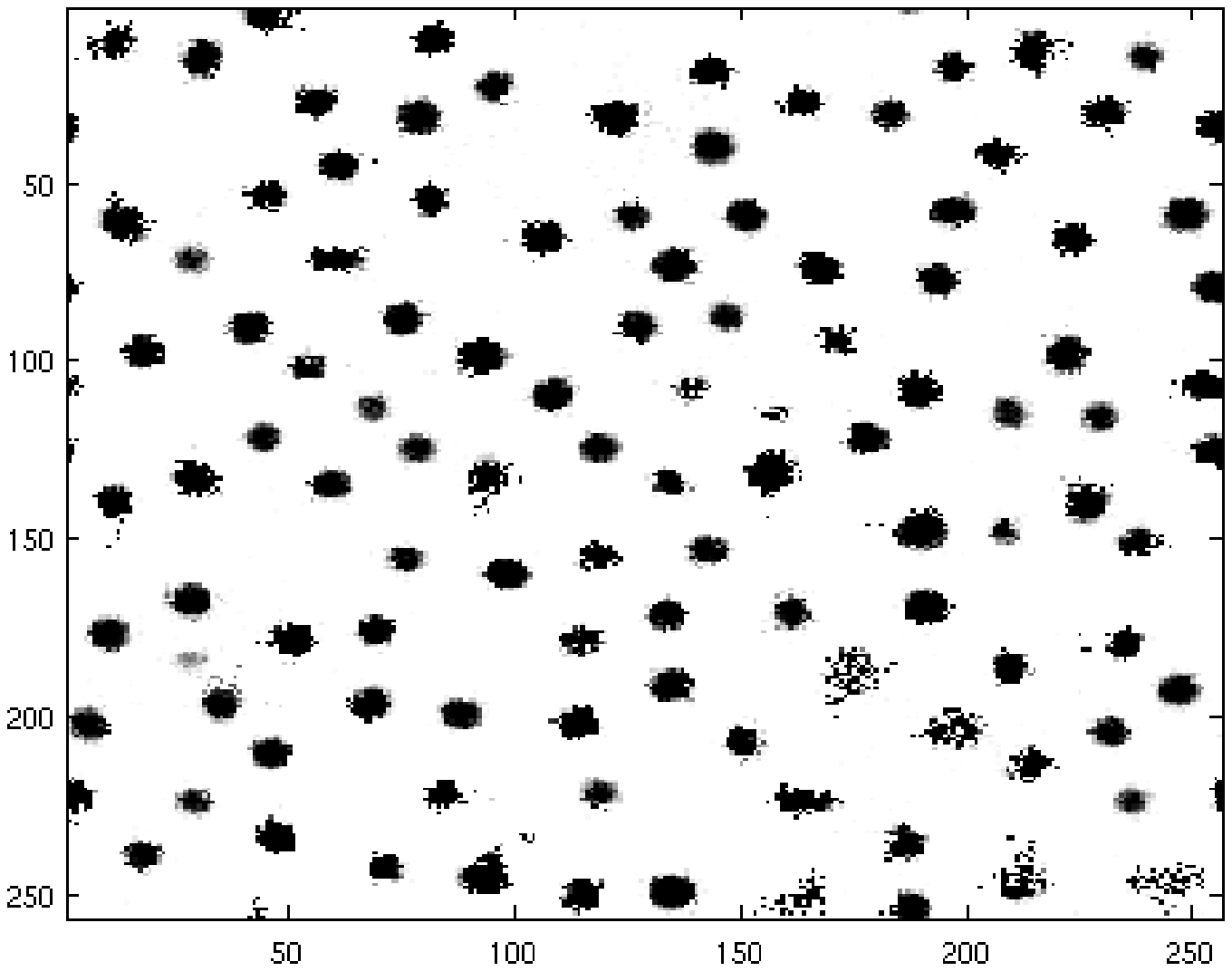}}{\unitlength=1bp}}\\
\subfloat[CGMC, q=4
]{{\label{fig:discs-cg4}}{\includegraphics[scale=0.25]{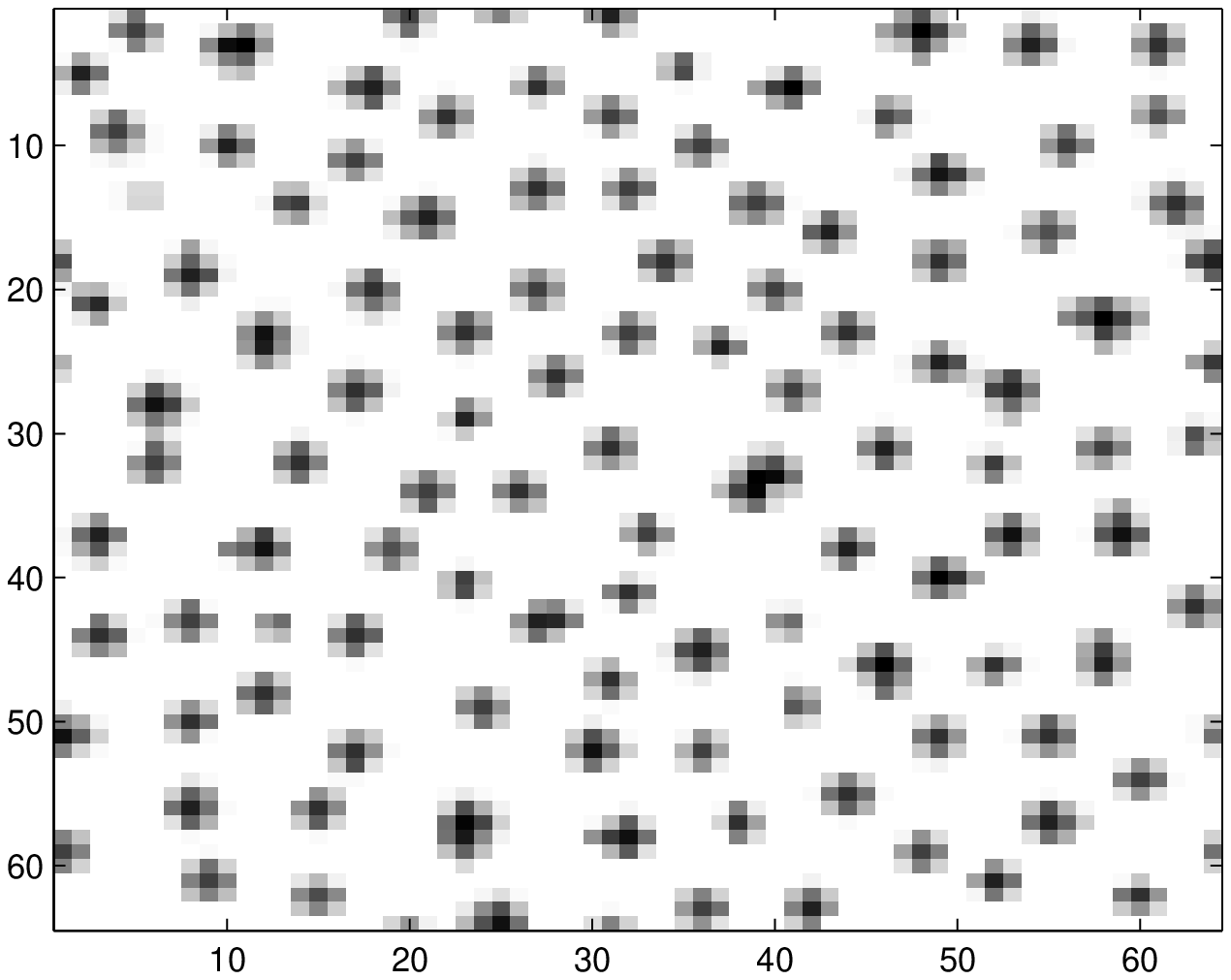}}{\unitlength=1bp}}
\subfloat[Two-level,
q=4]{{\label{fig:discs-co4}}{\includegraphics[scale=0.25]{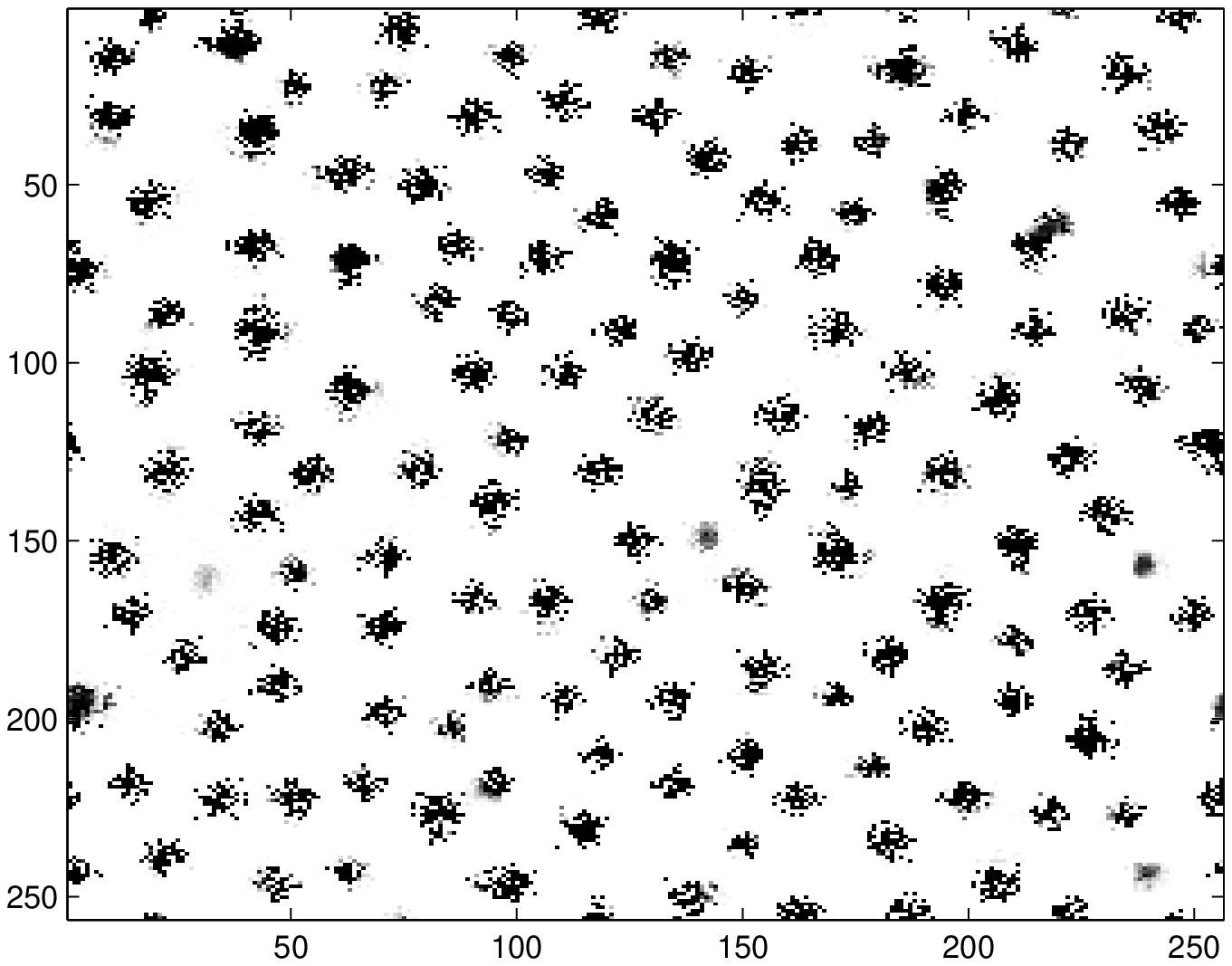}}{\unitlength=1bp}}
\caption{The two-level method with
correction is independent of the coarsening parameter.  Inverted discs $c_0=0.9 $, $\beta J_0= 0.6 $, $ N=256\!\times\! 256$. }
\label{fig:Discs}
\end{figure}

\begin{figure}[htp]
\centering
\subfloat[MH]{{\label{fig:stripes-tr}}{\includegraphics[scale=0.25]{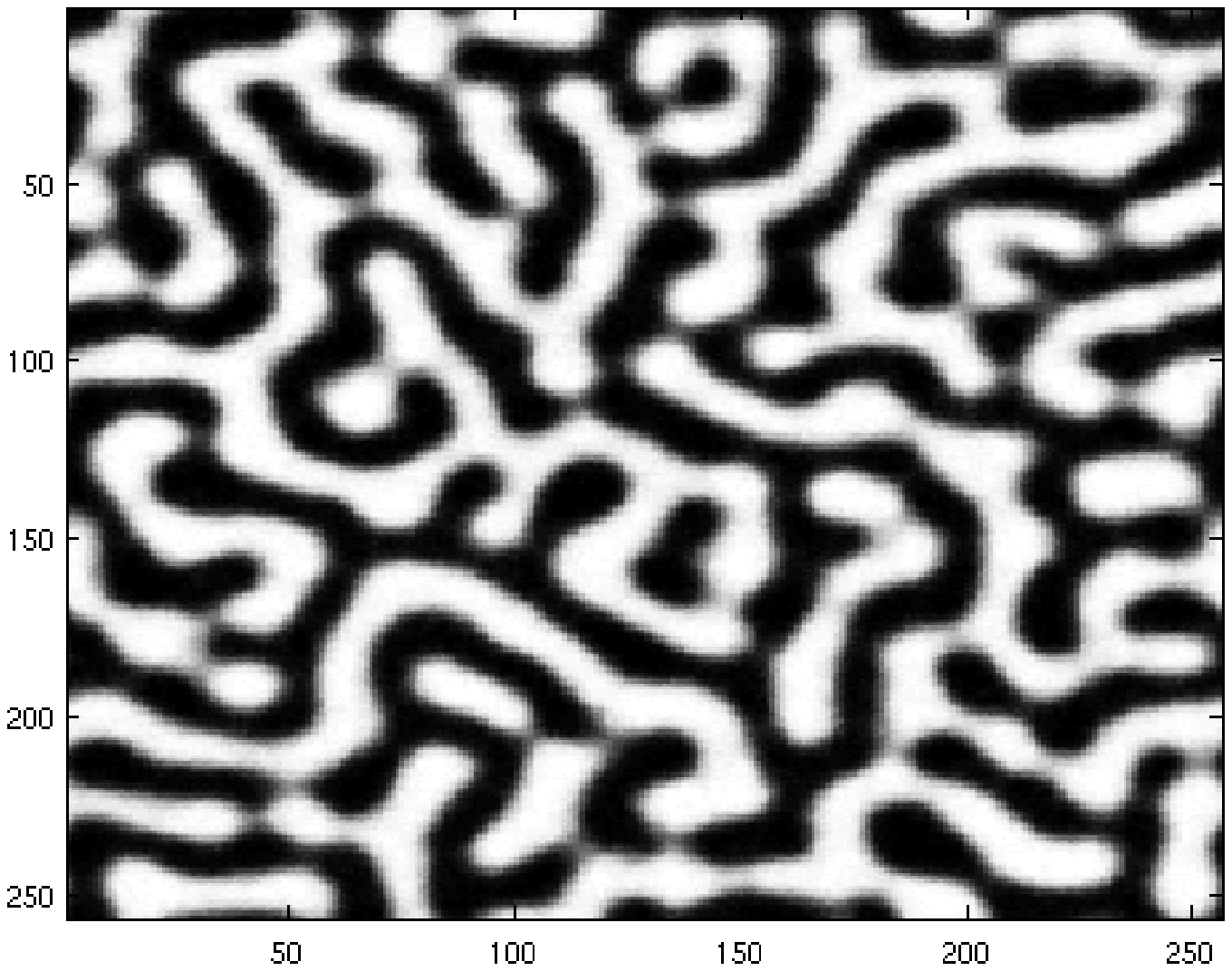}}{\unitlength=1bp}}
\subfloat[CGMC,
q=8]{{\label{fig:stripes-cg}}{\includegraphics[scale=0.25]{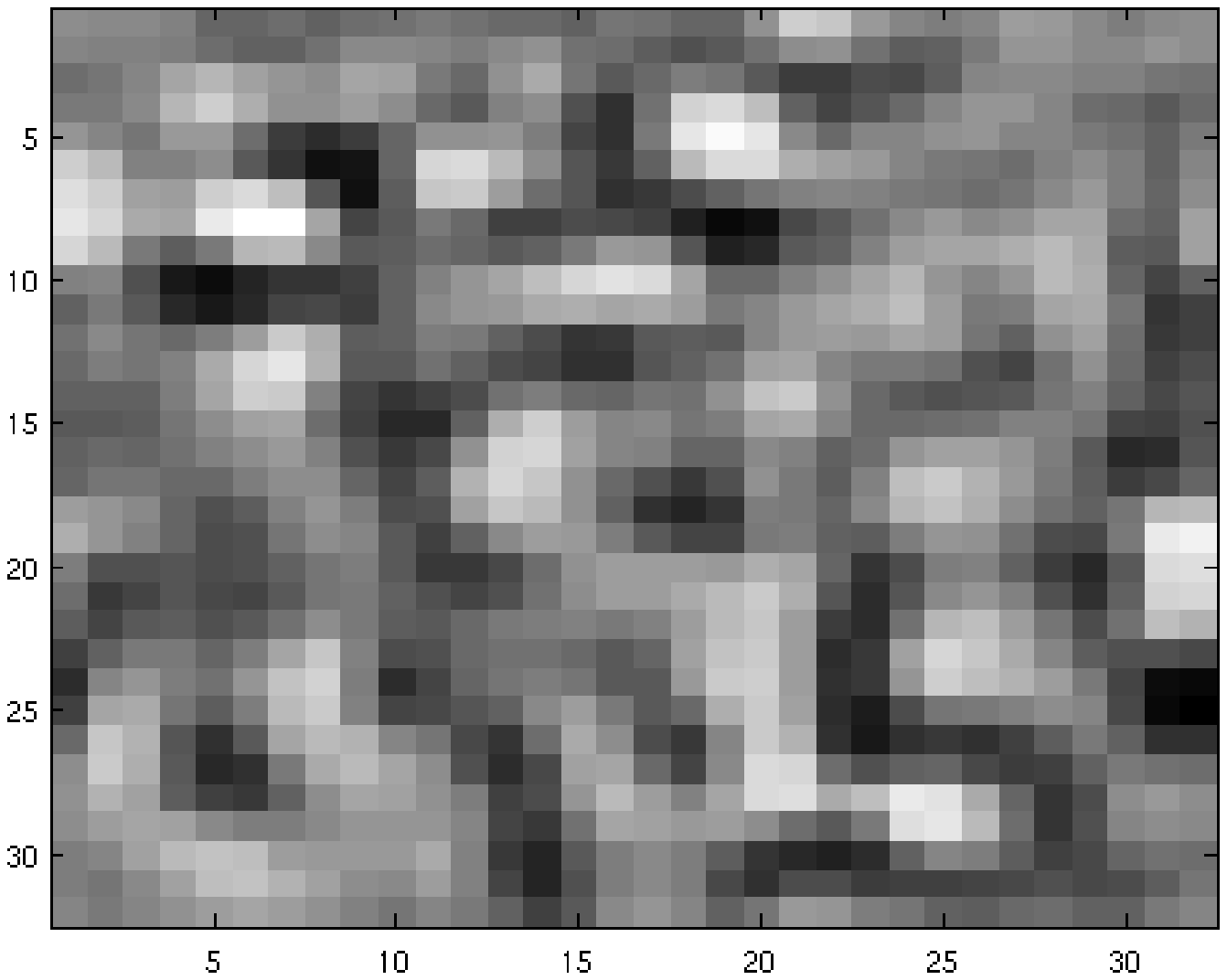}}{\unitlength=1bp}}
\subfloat[Two--level,
q=8]{{\label{fig:stripes-co}}{\includegraphics[scale=0.25]{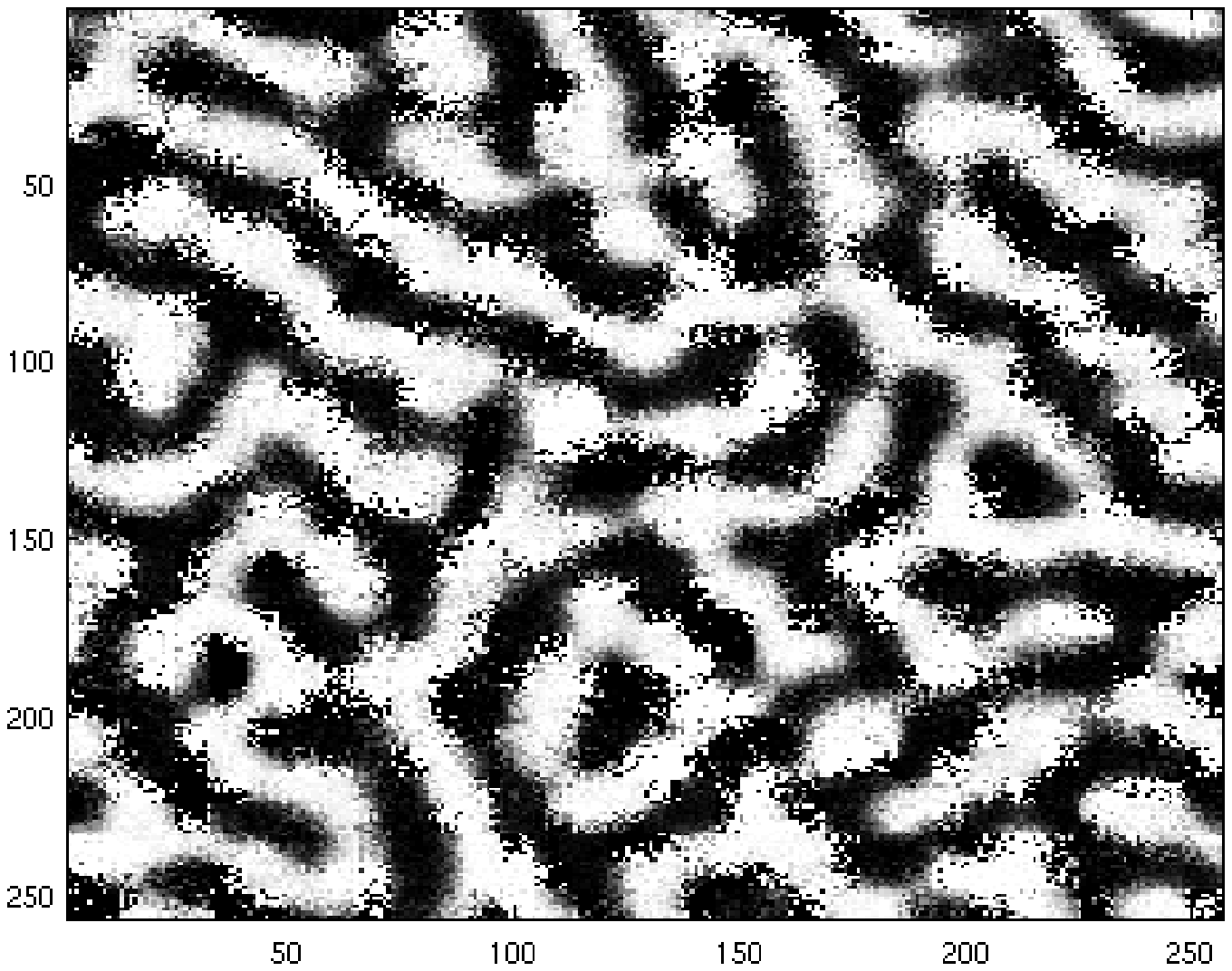}}{\unitlength=1bp}}\\
\subfloat[CGMC,
q=4]{{\label{fig:stripes-cg4}}{\includegraphics[scale=0.25]{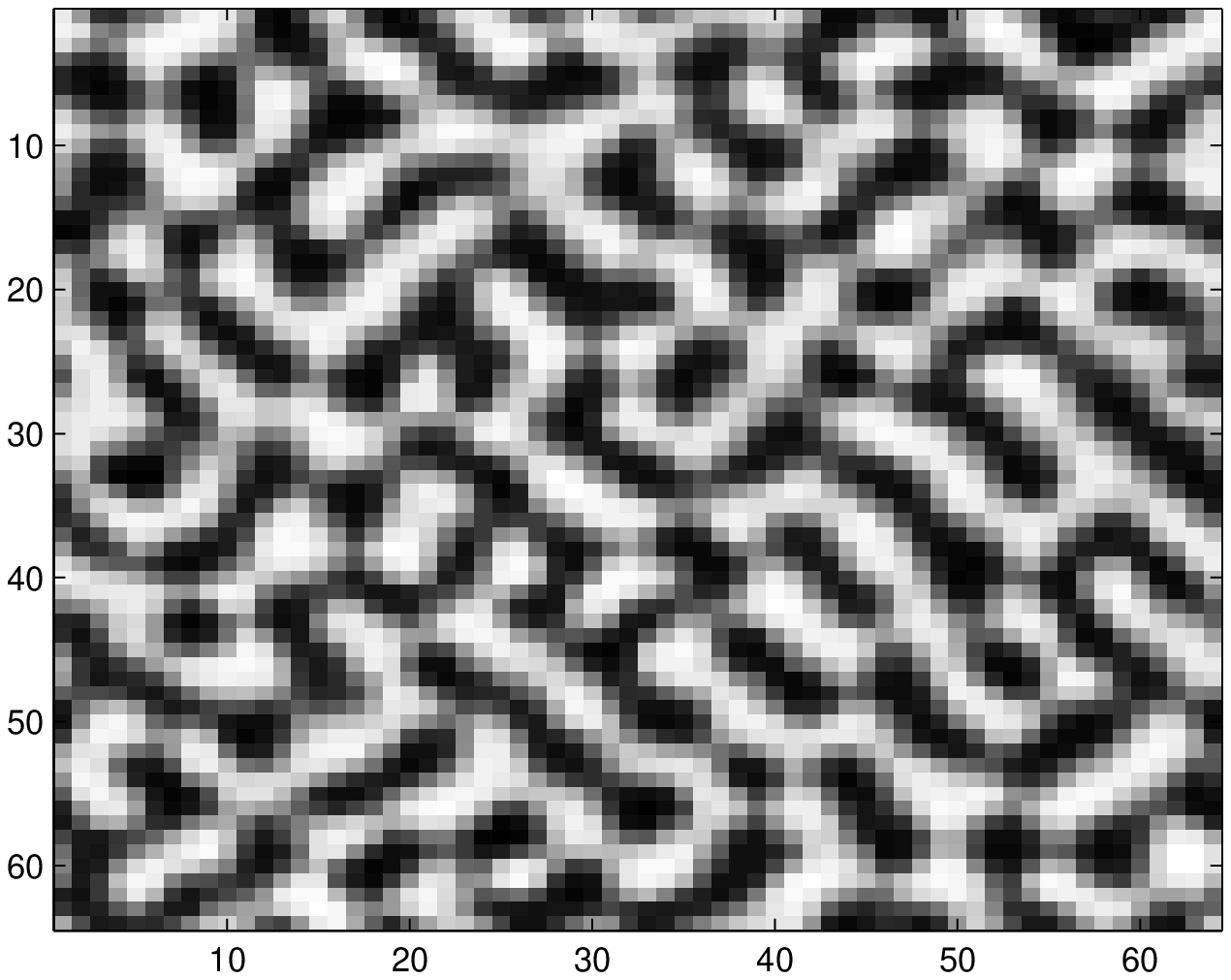}}{\unitlength=1bp}}
\subfloat[Two-level,
q=4]{{\label{fig:stripes-co4}}{\includegraphics[scale=0.25]{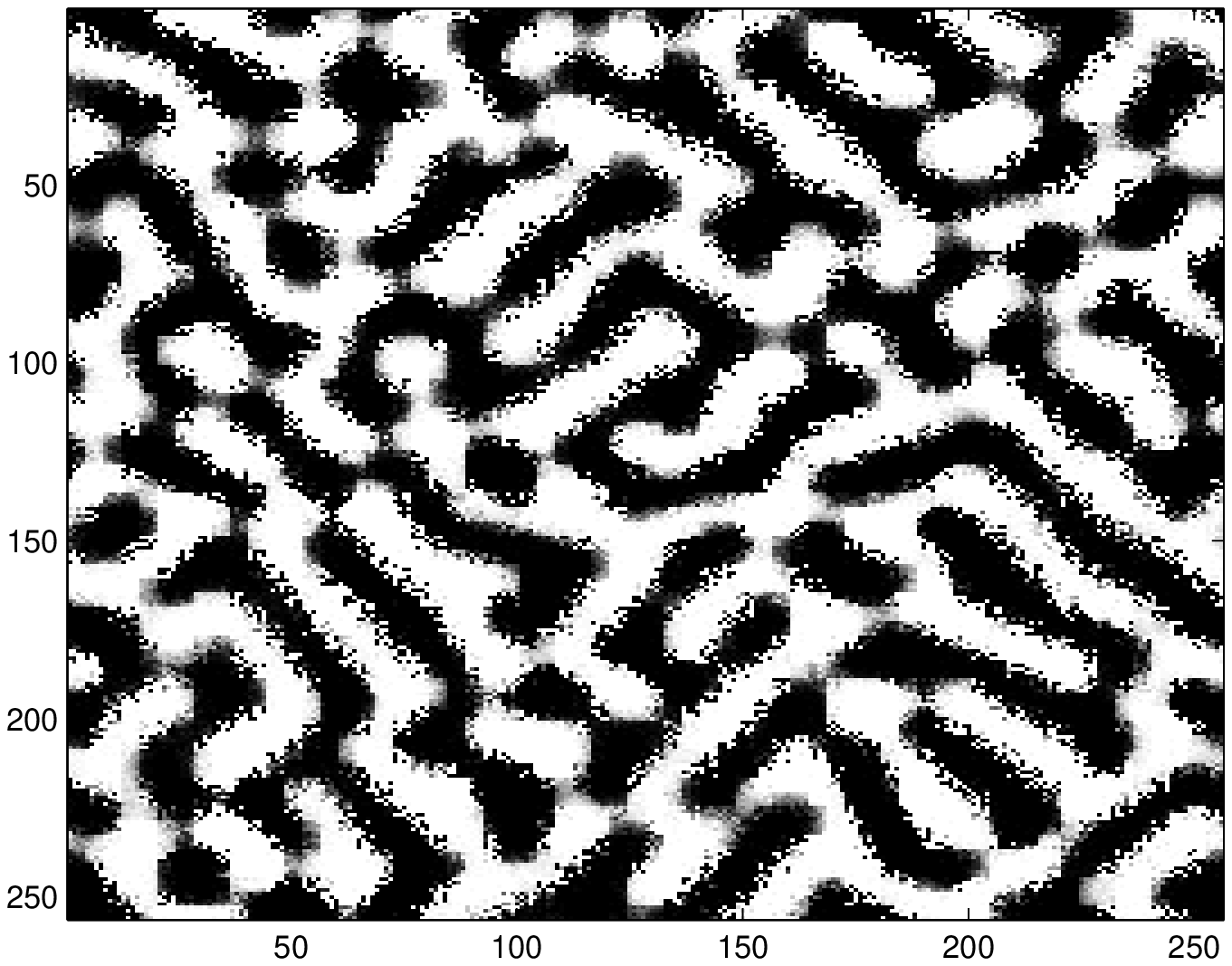}}{\unitlength=1bp}}
\caption{The two--level method   captures the
correct feature scaling even when CGMC is inaccurate. Labyrinths    $c_0=0.5 $, $\beta J_0= 0.2 $, $ N=256\!\times\! 256$.}
\label{fig:Stripes}
\end{figure}

In the tables and figures following we provide simulation results from
three methods, the conventional MH, the CGMC and the two--level CGMC.
For the latest method    two approaches were tested, the sampling with correction terms, Figures
 \ref{fig:Discs} - \ref{fig:Stripes} and Tables \ref{Diameter} - \ref{Methods}  and the approximate CG with potential splitting, Figure \ref{fig:Discs}. Figure \ref{fig:3} shows a graph of the interaction potential $ J(|x-y|) $   and  $J_c(|x-y|) =J(|x -y|)-\bar J (|x-y|)$  appearing in the second step of the two--level method,  representing the correction of  compressing at the first step of the method.    $\bar J(|x-y|) =\bar J(k,l) , \ x, y \in \LATT $ with $x\in C_k, y\in C_l $ is the compressed potential,  used for the coarse  space simulations.   
In all simulations presented in the sequel the range  of pure attractive and   repulsive  forces are $r_a=4.47$,  $r_r=10$ respectively  with repulsion strength $\chi=0.1$.  Since the interactions
follow an exponential decay, in  implementations  we  use  a cut--off range $L$ for the  potential $J(r)$ such that $J(r)< 10^{-6}$ for $ r> L $, which accounts on $ (2L)^2$ interactions for each site.   

For  total  coverage $ c_0=0.9$, that leads to  pattern  composed by periodic    inverted discs, Table \ref{Diameter} gives the average discs diameter and computational time for different  values of  lattice size $ N$,  obtained with the two--level CGMC method with the corrections strategy. The results   verify that the proposed method provides the expected behaviour of the system for a finite lattice  as its  size increases, that is  features properties stabilize. Features statistics are calculated with the use of edge detection techniques of image processing \cite{ParkerImPr}.  In order to demonstrate the effect of the statistical errors  we present  confidence interval estimates. Such intervals are obtained after transforming the simulations  data to follow an approximately normal distribution.  In simulations for all methods  the number of MC iterations  is $5\times10^7$.
\begin{table}[h]
\caption{Finite system  size behaviour  of the two--level CGMC method.
Average discs diameter $\AVG{d}$ and computational time. $c_0=0.9 $, $\beta J_0= 0.6 $ $q=8$, $L=24$. }
 \centering  
\begin{tabular}{l c c c c c  }  
\hline\noalign{\smallskip} 
&Lattice size  & Diameter  & Std & Confidence Interval & CPU
time(min) \\
 \noalign{\smallskip}\hline\noalign{\smallskip}
 &    $128\!\times\!128 $   &   10.2  &  3.6  &  [5.7,  13.3]   &   30  \\   
&     $256\!\times\! 256$   &   9.5   &  2.5  &   [8.7, 9.9]  &   43 \\
&     $512\!\times\! 512$   &   8.8   &  3.4  &   [8.5, 9.1]  &   66  \\  
 &    $1024\!\times\! 1024 $ &   8.6   &  3.6  & [8.5, 8.7] &   110 \\  
\\
\hline\noalign{\smallskip}   
 &    MH method      &        &   &     \\  
 \noalign{\smallskip}\hline\noalign{\smallskip} 
 &    $512\!\times\!512$ &   8.27  &  1.2  & [8.1, 8.3] &    805  \\  
\end{tabular}
\label{Diameter}
\end{table} 
For the sake of comparison in Table \ref{Diameter} we present also results obtained 
with the   MH method  for  a $N=512\!\times\!512$ lattice. For the same lattice the CGMC  method  was tested with coarse graining parameter $q=8$ that shows  significant  deviations,
 the average diameter   being  $\AVG{d}\sim q \!\times\! \AVG{\bar d} = 8\!\times\! 7.8 $.  This proves  that we {\it over-coarsened},  on the other hand for smaller values of the parameter $Q$ CGMC method provides better approximations,
for example when $q=2 $  $\AVG{\bar d}= 4.9$ and $\AVG{d}\sim 2\!\times\! 4.9$,  but the computational time, being $187 $ minutes, is still relatively large.

  Figures \ref{fig:Discs} and  \ref{fig:Stripes} show the averaged equilibrium conformations for  two total coverage values $ c_0=0.9$ and $c_0=0.5$, leading to the formation of inverted discs and stripes respectively, for  a $ N=256 \!\times\! 256 $ lattice.  Black and white  dots indicate occupied   and vacant sites respectively.     CGMC sampling is fast (i.e. $ 7.3$ minutes for $ q=8$) but with large   deviations from the expected pattern are significant, while the two--level CGMC refines that error providing correct feature scales, of course  with the cost of  increasing   the computational, see  Figures \ref{fig:stripes-cg} and \ref{fig:stripes-co}.  On the other hand there is a substantial reduce on computational time  compared to  the MH method,  from  $365$  to  $25 $ minutes.
For  a smaller coarse graining parameter,  $ q =4$ in  Figure \ref{fig:stripes-cg4},  the CGMC prediction is better   albeit at a higher computational cost where again and the coarsening error is   improved by the proposed method and microscopic details are available.

Table \ref{Methods} shows a comparison of the discs statistics and
computational  time for the conventional and the two--level method for various
coarse graining parameters.  The two--level method estimates are close to the MH even
for high values of the  coarse parameter, with a significant reduction of the
computational time.
 
\begin{figure}[htp]
\centering
\subfloat[MH]{{\label{fig:stripes-trs}}{\includegraphics[scale=0.25]{figures/TraditionalStripes}}{\unitlength=2bp}}
\subfloat[Splitting]{{\label{fig:stripes-spl}}{\includegraphics[scale=0.25]{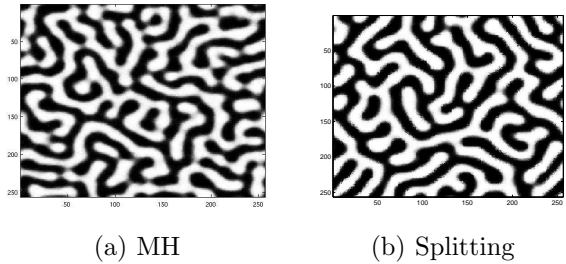}}{\unitlength=1bp}}
\caption{ Labyrinths. Coupled method with splitting potential, $ q=8,
s=4$. Splitting captures the qualitative picture but misses the
characteristic scale. }
\label{fig:StripesSplit}
\end{figure}

\begin{table}[htp]
\caption{ Two--level CGMC with correction, $N=256\!\times\! 256$} 
\centering  
\begin{tabular}{l c c c c c c }  
\hline\noalign{\smallskip} 
 & Method		   & Diameter $\AVG{d}$  &  Std    &CPU time (min) \\  
 \noalign{\smallskip}\hline\noalign{\smallskip}
& MH 		     &  8.7            &  1.3         &    252     \\   
& Two--level   q=2 &  8.4            &  1.6         &    66         \\
& Two--level   q=4 &  7.9            &  4.6         &    22       \\
& Two--level   q=8 &  8.9            &  3.2         &    23     \\
\hline  
\end{tabular}
\label{Methods}
\end{table} 
  Figures \ref{fig:StripesSplit} present results of simulations with the  approximate splitting approach  where  we  split the interactions up to  range $S$ and neglect the correction terms at the second step, see Section~\ref{Decompositions}. This strategy captures the qualitative picture but misses the characteristic length. From the simulations we tested we conclude that this approach  is not suggested  since the splitting disturbs significuntly the relative strength of attractive and repulsive forces that  has an impact on  characteristic length scales of the features.   

\section{Conclusions}
In this work we propose the multilevel CGMC method and study its
efficiency theoretically and numerically. The hierarchy of CGMC methods introducing multilevel
decompositions are shown to reduce the computational complexity of conventional
methods. The {\it reconstructing} measures, implemented by an 
accept/reject method, provide the correct  inverse mapping  from  coarse  to fine resolutions refining CGMC sampling.
With Theorem~\ref{main} we prove that in terms of speed of convergence
the two--level CGMC   is comparable to a conventional MH, that is
controlled by the coarse approximating measure $ \BARMZ(d\eta)$. In general
the method's mixing time cannot be better than a conventional method,
but the reduction of computational effort with the error control of the
coarse--graining approximation make it an adequate method for studying high
dimensional systems.
We successfully apply the method in a nano-pattern discovery problem,
verifying that the multilevel CGMC method provides the expected
behaviour of the system and microscopic information, even when CGMC is inaccurate.
 
\section*{Acknowledgements}


\appendix
\section{Proofs}\label{Math}
 Before studying the proposed method's mathematical properties we need to introduce some
definitions and theoretical facts \cite{RC}. Let $\{X_n\}$ be a Markov chain on
space $\Sigma$ with transition kernel $\Kk$.
\begin{defin}
\begin{enumerate}[i)]
\item A transition kernel $\Kk$ has the {\it stationary measure} $\mu$ if  
\[ \int_{\Sigma}\Kk(\sigma,\sigma') \mu(d\sigma')=\mu(\sigma), \ \ \text{for
all } \sigma\in \Sigma\PERIOD\]
\item $\Kk$ is called reversible with respect to $\mu$ if 
\[ (g,\Kk f)_{\mu}=(\Kk g, f)_{\mu},\ \ \text{ for all } g, \ f\in
L^2(\mu)\PERIOD\]
 \end{enumerate}
 \end{defin}
where $ (g,f)_{\mu} =\int_{\Sigma}\overline{g(\sigma)}f(\sigma)\mu(d\sigma)$,
$\overline{g}$ denoting the complex conjugate of $g$ and
$\Kk g(\sigma)=\int_{\Sigma}\Kk(\sigma,d\sigma')g(\sigma'), \forall \sigma \in
\Sigma$.
  A sufficient condition for $\mu$ being a stationary measure of $\Kk$ that is often easy to check is the detailed balance (DB) condition.
 \begin{defin}\label{defDB}{\bf (Detailed Balance)}
A Markov chain with transition kernel $\Kk$ satisfies the detailed balance
condition if there exists a function $f$ satisfying
 \begin{equation}\label{dBD}
\Kk(\sigma,\sigma')f(\sigma)=\Kk(\sigma',\sigma)f(\sigma'),\ \ \text{for all }
\sigma, \sigma' \in \Sigma \PERIOD
 \end{equation}
 \end{defin}

 Here we continue with the proof of  Theorem~\ref{ReversibilityThm}.
\pf 
{\it i)} 
 Let $\sigma\neq\sigma' $, recalling  the definition of transition kernel $ \Kk_{CG}(\sigma,\sigma')$ \VIZ{KCG} we have
\begin{eqnarray*}
\Kk_{CG}(\sigma,\sigma')\mu(\sigma) &=&  
\alpha_f(\sigma,\sigma') \alpha_{CG}(\COP\sigma,\COP\sigma')\mu_r(\sigma'|\COP\sigma')\BARIT{\rho}(\COP\sigma,\COP\sigma')\mu(\sigma)\\
&=& \min\left\{1,
\frac{\mu(\sigma')\BARMZ(\COP\sigma)\mu_r(\sigma|\COP\sigma)}{\mu(\sigma)\BARMZ(\COP\sigma')\mu_r(\sigma'|\COP\sigma')}\right\}\times\\
&& \!\!\!\!\!\min\left\{1, \frac{\BARMZ(\COP\sigma')\BARIT{\rho}(\COP\sigma',\COP\sigma)}{\BARMZ(\COP\sigma)\BARIT{\rho}(\COP\sigma,\COP\sigma')}\right\}  \mu_r(\sigma'|\COP\sigma') \BARIT{\rho}(\COP\sigma,\COP\sigma')\mu(\sigma)\\
&=&  \min\left\{1,
\frac{\mu(\sigma)\BARMZ(\COP\sigma')\mu_r(\sigma'|\COP\sigma')}{\mu(\sigma')\BARMZ(\COP\sigma)\mu_r(\sigma|\COP\sigma)}\right\} \times\\
&&\!\!\!\!\!\min\left\{1, \frac{\BARMZ(\COP\sigma)\BARIT{\rho}(\COP\sigma,\COP\sigma')}{\BARMZ(\COP\sigma')\BARIT{\rho}(\COP\sigma',\COP\sigma)}\right\} \mu_r(\sigma|\COP\sigma)\BARIT{\rho}(\COP\sigma',\COP\sigma)  \mu(\sigma')\\
&=& \Kk_{CG}(\sigma',\sigma)\mu(\sigma') \PERIOD
\end{eqnarray*}

{\it ii)} follows from   {\it i)}. Detailed balance with $ \mu(\sigma)$ is sufficient to guarantee that   $ \mu(\sigma)$  is the stationary distribution of kernel $ \Kk_{CG}(\sigma,\sigma')$, (Theorem 6.46 in \cite{RC}).

\noindent
 {\it iii)}  
To prove that chain $\{Y_n\}$ is $\mu$-irreducible we need to prove
$\Kk_{CG}(\sigma,A) > 0$, for all $\sigma \in E $ and $A$ measurable such that
$ \mu(A)>0$.  We have
 \begin{align*}
\Kk_{CG}(\sigma,A) &= \int_A \Kk_{CG}(\sigma,\sigma') d\sigma'\ge
\int_{A-\{\sigma\}} \Kk_{CG}(\sigma,\sigma') d\sigma' =\\
&= \int_{A-\{\sigma\}} \alpha_f(\sigma,\sigma') \alpha_{CG}(\COP
\sigma, \COP \sigma')
\mu_r(\sigma'|\COP\sigma')\bar{\rho}(\COP\sigma,\COP\sigma')
d\sigma'\PERIOD
 \end{align*}
From assumptions on $\bar \rho(\eta,\eta') $ and $\mu_r(\sigma|\eta) $ term $
\mu_r(\sigma'|\COP\sigma')\bar{\rho}(\COP\sigma,\COP\sigma')$
is positive for all $ \sigma, \sigma' \in E$. Also since $ A \subset E $ and $E \subset supp(\mu_0) $,
$ \alpha_f(\sigma,\sigma') $ and $ \alpha_{CG}(\COP\sigma,\COP\sigma')
$ are positive. These ensure that $ \Kk_{CG}(\sigma,A) >0$.
 
 \noindent
{\it iv)} A sufficient condition ensuring that $ \{ Y_n\}$ is aperiodic is that
$\Kk(\sigma,\{\sigma\}) >0 $ for some $ \sigma \in E$, that means the event $
Y_{n+1}=Y_n$ happens with positive probability. We have
 \begin{align*}
\Kk_{CG}(\sigma, \{ \sigma \} ) = 1 - \int_{ \{\sigma' \neq \sigma\}}
\alpha_f(\sigma,\sigma') \alpha_{CG}(\COP \sigma, \COP \sigma')
\mu_r(\sigma'|\COP\sigma')\bar{\rho}(\COP\sigma,\COP\sigma')
d\sigma'\PERIOD
\end{align*}
If for all $\sigma\in \Sigma $, $ \Kk_{CG}(\sigma, \{ \sigma \} ) = 0$
then
\begin{equation*}
\int_{ \{\sigma' \neq \sigma\}} \alpha_f(\sigma,\sigma')
\alpha_{CG}(\COP \sigma, \COP \sigma')
\mu_r(\sigma'|\COP\sigma')\bar{\rho}(\COP\sigma,\COP\sigma') d\sigma'
=1\COMMA
 \end{equation*}
which means that $\alpha_f(\sigma,\sigma') =1$ and $\alpha_{CG}(\COP
\sigma, \COP \sigma') =1 $ for almost all $ \sigma \in \{ \sigma \in E
: \bar\rho(\COP\sigma,\COP\sigma')\mu_r(\sigma'|\COP\sigma') > 0,\
\text{for some } \sigma'\in E\}$.
This would mean that the reconstructed proposal  kernel
$\bar\rho(\COP\sigma,\COP\sigma')\mu_r(\sigma'|\COP\sigma') $ is
sampling from the exact target measure $ \mu$ which in general is not
true. Therefore there exists $\sigma \in E $ such that $
\Kk_{CG}(\sigma, \{ \sigma \} ) = 0$.
   \hfill $\square$

 \section{Benchmark example calculations }\label{appendix1}
The detailed calculations of the application of
Theorem~\ref{main} for the benchmark example described in Section
\ref{benchmark} are provided in detail here.

The conventional MH algorithm is described by a  transition  kernel  proposing a spin-flip at the site x with the probability $1/N$,  that is $\rho(\sigma,\sigma') =N^{-1}\sum_{x\in\LATT}\delta(\sigma'-\sigma^x)$,
 with  acceptance probability
\begin{equation*}
\alpha(\sigma,\sigma^x)=\min\left\{1, e^{(1-2\sigma(x))[K
\left(\sigma(x-1)+\sigma(x+1)\right)+\frac{J}{N}\sum_{y\neq x}^{
N}\sigma(y)]}\right\} \COMMA
\end{equation*}
where, for simplicity,  we consider the case of zero external field $h$.
In the   two--level CGMC   algorithm the proposal probability distribution is
\begin{equation*}
\bar{\rho}(\eta,\eta') = \frac{1}{M}\sum
_{k\in\LATTC}\left[\frac{Q-\eta(k)}{Q} \delta(\eta'-\eta^{k+}) +
\frac{\eta(k)}{Q}\delta(\eta'-\eta^{k-})\right]\COMMA\]
where $\eta^{k\pm}(l)= \eta(l), l\neq k$ and $\eta^{k\pm}
(k)=\eta(k)\pm 1$ and $k\pm $ denotes the adsorption $(+)$ or desorption $(-)$
in the cell $k$. The acceptance probability of the first level is
\begin{equation*}
\alpha_{CG}(\eta,\eta^{k\pm})=\min\left\{1, e^{- 
\Delta_{k\pm}\BARH_l(\eta) } \right\}  \PERIOD
\end{equation*}
where $ \Delta_{k+}\BARH_l(\eta) = -\frac{J}{N} \sum_{l\in\LATTC} \eta(l)   $
and $ \Delta_{k-}\BARH_l(\eta) = -\frac{J}{N}\left[1- \sum_{l\in\LATTC}
\eta(l) \right] $.
The reconstruction probability distribution is a uniform distribution in each
cell $\mu_r(\sigma^x|\eta^{k\pm}) = \prod_{l\in \LATTC}
\mu^{(l)}_r(\sigma^x_{C_l}| \eta^{k\pm}(l))$. Note that at each MC iteration
the change in the new state $ \eta^k$ happens in cell $k$, thus we need to perform the
reconstruction of $ \sigma^x$ only in the cell $k$, according to
\begin{equation*}
\mu_r(\sigma^x_{C_k}|\eta^k(k)) =
\frac{1 }{Q-\eta(k)}\delta(k-k+) +
\frac{1}{\eta(k)}\delta(k-k-) \PERIOD
\end{equation*}
According to Algorithm \ref{mCGMC}, the second level acceptance   probability is
\begin{equation*}
 \alpha_{f}(\sigma,\sigma^x)=\min\left\{1, e^{-  \Delta_x
H_s(\sigma) }\right\} =\min\left\{1, e^{ K (1-2\sigma(x))
\left(\sigma(x-1)+\sigma(x+1)\right)  } \right\} \PERIOD
\end{equation*}
Term $ \mathcal{B}(\sigma,\sigma')$ \VIZ{DefB} is   equal to one for all $ \sigma,\sigma'\in \Sigma_N$. Indeed for $\eta'=\eta^{k+}$ and $ x\in C_k \subset \LATT$,  
\begin{equation*}
\frac{\bar{\rho}(\eta,\eta^{k+})
\mu_r(\sigma^x|\eta^{k+})}{\rho(\sigma,\sigma^x)} =
\frac{\frac1M \frac{Q-\eta(k)}{Q}  \frac{1}{Q-\eta(k)} }{\frac1N}=1\PERIOD
\end{equation*}
Similarly for $\eta'=\eta^{k-}$,
\begin{equation*}
\frac{\bar{\rho}(\eta,\eta^{k-})
\mu_r(\sigma^x|\eta^{k-})}{\rho(\sigma,\sigma^x)} =
 \frac{\frac1M \frac{\eta(k)}{Q}  \frac{1}{\eta(k)} }{\frac1N}=1
\end{equation*}
Therefore its upper and lower bounds are $\underline{\gamma}=\bar{\gamma}=1$.
We consider the splitting approach, and use that $H_N(\sigma) -\BARH_l(\eta) =
H_s(\sigma)$ and $H_l(\sigma) = \BARH_l(\eta) $. The for any $x\in\LATT$ such
that $ x\in C_k,\ k\in \LATTC$ we have
\begin{equation*}
\mathcal{A}(\sigma,\sigma^x)=
\left\{  \begin{array}{cc}
    1,\ \  &\text{ if } (\sigma,\sigma^x)\in  C_1 \\
\min\{ e^{-\beta\Delta_k\BARH_l(\eta)},
e^{\beta\Delta_k\BARH_l(\eta)} \} , \ \ &\text{ if }
(\sigma,\sigma^x)\in C_2 \\
\min\{ e^{-\beta\Delta_x H_s(\sigma)}, e^{\beta\Delta_x H_s(\sigma)} \} , \ \
&\text{
if } (\sigma,\sigma^x)\in C_3
      \\
\min\{e^{-\beta\Delta_x H_N(\sigma)}, e^{\beta\Delta_x H_N(\sigma)} \} , \ \
&\text{ if }
(\sigma,\sigma^x)\in C_4
\end{array}  \right.
\end{equation*}
 Set $C_4=\emptyset$ according to the following argument.
Let $ (\sigma, \sigma^x)  \in C_4$ such that $\alpha_{CG}=1,\ \alpha_f=1$ and $
\alpha<1$.The first two relations are equivalent to $\Delta_k\BARH_l(\eta)\le 0$
and $\Delta_x H_s(\sigma)\le 0$ that imply $\Delta_x H_N(\sigma) \le 0 $, since
$ H_N(\sigma) = H_s(\sigma) + \BARH_l(\eta)$, thus $\alpha =1 $, a
contradiction. Analogous argument holds for the case $ \alpha=1, \alpha_{CG}<1,
\alpha_f<1 $ that proves $C_4=\emptyset$. 

Let us consider  $ K\neq 0$. Using the analytic expression of $\Delta_x
H_s(\sigma) $ and $\Delta_k \BARH_l(\eta) $  
we have, for $(\sigma,\sigma^x) \in C_2 $
\begin{eqnarray*}
\min\left\{ e^{-\beta\Delta_{k+}\BARH_l(\eta)},
e^{\beta\Delta_{k+}\BARH_l(\eta)} \right\}& \ge& e^{-| \beta
\Delta_{k+}\BARH_l(\eta)|}|_{\{\eta(l)=Q, \forall l\in \LATTC\}} \\
&=& e^{-|\frac{J}{N} N|}=e^{-|J|}\PERIOD
\end{eqnarray*}
Similarly for $(\sigma,\sigma^x) \in C_3 $
\begin{eqnarray*}
&\min\left\{ \EXP{-\beta\Delta_{k-}\BARH_l(\eta)},
\EXP{\beta\Delta_{k+}\BARH_l(\eta)} \right\} \ge \EXP{-|\frac{J}{N}
(N-1)|}>\EXP{-|J|} \text{ and }\\
&\min\left\{ \EXP{-\beta\Delta_{x}\BARH_s(\sigma)},
\EXP{\beta\Delta_{x}\BARH_s(\sigma)} \right\} \ge  \EXP{-2|K|}\PERIOD
\end{eqnarray*}
Therefore
$
\inf_{\sigma,
x}\Aa(\sigma,\sigma^x)=\Aa(\sigma,\sigma^x)|_{\{\sigma:
\sigma(x)=1, \forall x \in \LATT\}}=\min\{\EXP{-|J|},\EXP{-2|K|} \}\PERIOD
$

When $ K=0$  then  $C_2= \emptyset$,   that can be proved with simple arguments,
and  for $(\sigma,\sigma^x) \in C_3 $ $ \Aa(\sigma,\sigma^x)=1$ since   $ \Delta_x H_s(\sigma) = 0$. 
Proving that $ \Aa(\sigma,\sigma^x) =1, $ for all $x\in \LATT $ and $\sigma\in \Sigma_N$ . 
 



\bibliographystyle{model1b-num-names}

\end{document}